\documentclass[11pt,a4paper]{article}

\usepackage{url,amsmath,amssymb,latexsym,pstricks,mathrsfs,comment,amsthm,graphicx,tikz,tikz-cd,enumerate,accents,pgffor,cite,wrapfig,multicol,float,cases,enumitem,bibspacing,calc}

\usepackage[colorlinks]{hyperref}

\usepackage{geometry} 

\usepackage{hhline}
\usepackage{multirow}

\usepackage[T1]{fontenc}
\usepackage{ wasysym }
\usepackage{stmaryrd}

\usepackage{fnpct}
\setfnpct{dont-mess-around}

\newcommand\nc\newcommand
\nc\rnc\renewcommand

\newcounter{ncols}
\newcounter{incols}
\newenvironment{partn}[1]{
  \setcounter{ncols}{#1} \setcounter{incols}{\thencols - 1}\setlength{\arraycolsep}{1pt}
  \Bigl( \hspace{-1.5truemm}\scriptsize 
    \begin{array}{@{\hskip 3pt} c *{\theincols}{|c} @{\hskip 3pt}  }
}{
     \end{array}
     \normalsize \hspace{-1.5truemm}\Bigr)\setlength{\arraycolsep}{6pt}
}

\nc\ben{\begin{enumerate}[label=\textup{(\roman*)},leftmargin=7mm]}
\nc\bena{\begin{enumerate}[label=\textup{(\alph*)},leftmargin=7mm]}
\nc\BEN{\begin{enumerate}[label=\textup{(\Roman*)},leftmargin=7mm]}
\nc\een{\end{enumerate}}
\nc\bit{\begin{itemize}}
\nc\eit{\end{itemize}}
\nc\pf{\begin{proof}}
\nc\epf{\end{proof}}
\rnc\iff{\ \Leftrightarrow\ }
\rnc\implies{\ \Rightarrow\ }
\nc\an{\mathrel\&}
\nc\epfres{\hfill\qed}
\nc\ol\overline
\nc\B{\mathcal B}
\nc\PT{\mathcal P\mathcal T}
\nc\restr{{\restriction}}
\nc\De\Delta
\nc\id{\operatorname{id}}
\nc\Sub{\operatorname{Sub}}
\nc\Eq{\operatorname{Eq}}
\nc\sm\setminus
\nc\supp{\operatorname{supp}}
\nc\cosupp{\operatorname{cosupp}}
\nc\Supp{\operatorname{Supp}}
\nc\Cosupp{\operatorname{Cosupp}}
\nc\bmc{\begin{multicols}}
\nc\emc{\end{multicols}}
\nc{\pfitem}[1]{\medskip\noindent #1.}
\nc{\firstpfitem}[1]{#1.}
\nc\aftercases{\medskip\noindent}
\nc\Reg{\operatorname{Reg}}
\nc\Rad{\operatorname{Rad}}
\nc\Part{\operatorname{Part}}
\nc\Rest{\operatorname{Rest}}
\nc\PB{\P\B}
\nc\bA{{\bf A}}
\nc\bB{{\bf B}}
\nc\bC{{\bf C}}
\nc\bS{{\bf S}}
\nc\bX{{\bf X}}
\nc\F{\mathcal F}
\nc\N{\mathbb N}
\nc\K{\mathbb K}
\nc\C{\mathscr C}
\nc\Om\Omega
\nc\Si\Sigma
\nc\ze\zeta
\nc\nab\nabla
\nc\RR{\mathcal R}
\nc\LL{\mathcal L}
\nc\Pfd{\P^{\operatorname{fd}}}
\nc\Pfk{\P^{\operatorname{fk}}}
\nc\Pfcd{\P^{\operatorname{fcd}}}
\nc\Ffd{\F^{\operatorname{fd}}}
\nc\RP{\RR\P}
\nc\RJ{\RR\J}

\nc{\uv}[1]{\fill (#1,2)circle(.17);}
\nc{\lv}[1]{\fill (#1,0)circle(.17);}
\nc{\uvs}[1]{{\foreach \x in {#1} { \uv{\x}}}}
\nc{\lvs}[1]{{\foreach \x in {#1} { \lv{\x}}}}
\nc\colrect[5]{\fill[#5!20](#1,2)--(#2,2)--(#4,0)--(#3,0);}
\nc\colrecthigh[3]{\fill[#3!20](#1,1.75)--(#2,1.75)--(#2,2)--(#1,2);}
\nc\colrectlow[3]{\fill[#3!20](#1,.25)--(#2,.25)--(#2,0)--(#1,0);}
\nc{\uvw}[1]{\draw[fill=white] (#1,2)circle(.18);}
\nc{\lvw}[1]{\draw (#1,0)circle(.18);}
\nc{\uvws}[1]{\foreach \x in {#1}{ \uvw{\x}}}
\nc{\lvws}[1]{\foreach \x in {#1}{ \lvw{\x}}}
\nc{\uverts}[1]{\foreach \x in {#1}{ \uvert{\x}}}
\nc{\lverts}[1]{\foreach \x in {#1}{ \lvert{\x}}}
\nc{\uarcs}[1]{
{\foreach \x/\y in {#1}
{ \uarc{\x}{\y} }
}
}

\rnc\emptyset\varnothing

\nc{\darcs}[1]{
{\foreach \x/\y in {#1}
{ \darc{\x}{\y} }
}
}
\nc{\stlines}[1]{
{\foreach \x/\y in {#1}
{ \stline{\x}{\y} }
}
}

\nc\mt\mapsto
\nc\sub\subseteq
\nc{\COMMA}{,\qquad}
\nc\AND{\qquad\text{and}\qquad}
\nc\ANd{\quad\text{and}\quad}
\rnc{\H}{\mathrel{\mathscr H}}
\rnc{\L}{\mathrel{\mathscr L}}
\nc{\R}{\mathrel{\mathscr R}}
\nc{\D}{\mathrel{\mathscr D}}
\nc{\JJ}{\mathrel{\mathscr J}}
\nc{\J}{\mathcal J}
\nc{\LE}{\mathrel{\mathscr L}^E}
\nc{\RE}{\mathrel{\mathscr R}^E}
\nc{\HE}{\mathrel{\mathscr H}^E}
\nc{\LF}{\mathrel{\mathscr L}^F}
\nc{\RF}{\mathrel{\mathscr R}^F}
\rnc{\HF}{\mathrel{\mathscr H}^F}
\nc{\RG}{\mathrel{\mathscr R}^G}
\nc\upL[1]{{}^\uparrow#1}
\nc\upR[1]{#1^\uparrow}
\rnc\P{\mathcal P}
\rnc\S{\mathcal S}
\nc\Pt{\mathcal P t}
\nc\Bin{\mathcal B}
\nc\Dif{\mathcal D}
\nc\I{\mathcal I}
\nc{\leqR}{\leq_{\R}}
\nc{\leqL}{\leq_{\L}}
\nc{\leqJ}{\leq_{\JJ}}
\nc{\leqH}{\leq_{\H}}
\nc{\leqK}{\leq_{\K}}
\nc{\geqR}{\geq_{\R}}
\nc{\geqL}{\geq_{\L}}
\nc{\geqJ}{\geq_{\JJ}}
\nc{\geqH}{\geq_{\H}}
\nc{\geqK}{\geq_{\K}}
\nc\T{\mathcal T}
\nc\lam\lambda
\nc\rank{\operatorname{rank}}
\nc\dom{\operatorname{dom}}
\nc\codom{\operatorname{codom}}
\nc\im{\operatorname{im}}
\nc\coker{\operatorname{coker}}
\nc\Ker{\operatorname{Ker}}
\nc\Coker{\operatorname{Coker}}
\nc\set[2]{\{#1:#2\}}
\nc\bigset[2]{\big\{#1:#2\big\}}
\nc{\lar}[1]{ \xrightarrow {\ #1\ }}
\nc\pre\preceq
\nc\la\langle
\nc\ra\rangle

\nc\al\alpha
\nc\be\beta
\nc\ga\gamma
\nc\de\delta
\nc\ve\varepsilon
\rnc\th\theta

\nc{\uvert}[1]{\fill (#1,2)circle(.2);}
\rnc{\lvert}[1]{\fill (#1,0)circle(.2);}
\nc{\darcx}[3]{\draw(#1,0)arc(180:90:#3) (#1+#3,#3)--(#2-#3,#3) (#2-#3,#3) arc(90:0:#3);}
\nc{\darc}[2]{\darcx{#1}{#2}{.4}}
\nc{\uarcx}[3]{\draw(#1,2)arc(180:270:#3) (#1+#3,2-#3)--(#2-#3,2-#3) (#2-#3,2-#3) arc(270:360:#3);}
\nc{\uarc}[2]{\uarcx{#1}{#2}{.4}}
\nc{\stline}[2]{\draw(#1,2)--(#2,0);}

\nc{\custpartn}[3]{{\lower1.4 ex\hbox{
\begin{tikzpicture}[scale=.3]
\foreach \x in {#1}
{ \uvert{\x}  }
\foreach \x in {#2}
{ \lvert{\x}  }
#3 \end{tikzpicture}
}}}

\let\oldproofname=\proofname
\renewcommand{\proofname}{\rm\bf{\oldproofname}}

\allowdisplaybreaks

\begin{document}

\numberwithin{equation}{section}

\newtheorem{thm}[equation]{Theorem}
\newtheorem{lemma}[equation]{Lemma}
\newtheorem{cor}[equation]{Corollary}
\newtheorem{prop}[equation]{Proposition}
\newtheorem{conj}[equation]{Conjecture}

\theoremstyle{definition}

\newtheorem{rem}[equation]{Remark}
\newtheorem{defn}[equation]{Definition}
\newtheorem{eg}[equation]{Example}
\newtheorem{ass}[equation]{Assumption}
\newtheorem{prob}[equation]{Problem}
\newtheorem{ques}[equation]{Question}

\title{\vspace{-1.7cm}Ehresmann theory and partition monoids}
\author{}
\date{}

\maketitle
~\vspace{-2.1cm}
\begin{center}
{\large 
James East%
\footnote{\label{footnote:JE}Centre for Research in Mathematics and Data Science, Western Sydney University, Australia. {\it Email:} {\tt J.East\,@\,WesternSydney.edu.au}.  Supported by ARC Future Fellowship FT190100632.}
and Robert D.~Gray%
\footnote{School of Mathematics, University of East Anglia, United Kingdom. {\it Email:} {\tt Robert.D.Gray\,@\,uea.ac.uk}. Supported by EPSRC grant EP/N033353/1.}%
}
\end{center}

\begin{abstract}
This article concerns Ehresmann structures in the partition monoid $\P_X$.  Since $\P_X$ contains the symmetric and dual symmetric inverse monoids on the same base set $X$, it naturally contains the semilattices of idempotents of both submonoids.  We show that one of these semilattices leads to an Ehresmann structure on $\P_X$ while the other does not.  We explore some consequences of this (structural/combinatorial and representation theoretic), and in particular characterise the largest left-, right- and two-sided restriction submonoids.  The new results are contrasted with known results concerning relation monoids, and a number of interesting dualities arise, stemming from the traditional philosophies of inverse semigroups as models of partial symmetries (Vagner and Preston) or block symmetries (FitzGerald and Leech): ``surjections between subsets'' for relations become ``injections between quotients'' for partitions.  We also consider some related diagram monoids, including rook partition monoids, and state several open problems.


\end{abstract}

\vspace{-.5cm}
\tableofcontents

\section{Introduction}\label{s:intro}

Recall that a semigroup $S$ is inverse if for each element $x\in S$ there exists a unique $a\in S$ satisfying $x=xax$ and $a=axa$; this element, the inverse of $x$, is typically denoted ${a=x^{-1}}$.  As explained in~\cite{Lawson1998}, inverse semigroups model partial symmetries of mathematical structures; for a dual approach via block symmetries, see \cite{FL1998}.  The unary operation of inversion, ${S\to S:x\mt x^{-1}}$, leads to an equational characterisation \cite{Schein1963} of inverse semigroups as the variety of unary semigroups satisfying the laws
\[
(x^{-1})^{-1} = x \COMMA x=xx^{-1}x \COMMA (xy)^{-1} = y^{-1}x^{-1} \COMMA (xx^{-1})(yy^{-1})=(yy^{-1})(xx^{-1}).
\]
Dropping the last of the above laws leads to the class of regular $*$-semigroups \cite{NS1978}, key examples of which are the diagram monoids of Brauer, Jones, Kauffman, Martin, Temperley and Lieb \cite{HR2005,Brauer1937,TL1971,Kauffman1990,Jones1994_2,Martin1994}.  Diagram monoids (and associated algebras and categories) arise in many fields of mathematics and science---including invariant theory, classical groups, representation theory, logic, knot theory and statistical mechanics---and are the main motivating examples for the current work.  For more background on these monoids we refer to our previous article \cite{EG2017}.

Arguably one of the most important results in inverse semigroup theory is the Ehresmann-Schein-Nambooripad Theorem \cite[Theorem 4.1.8]{Lawson1998}, which establishes an isomorphism between the categories of inverse semigroups and inductive groupoids.  This celebrated result was a major impetus for the definition of Ehresmann semigroups and categories, given by Lawson in \cite{Lawson1991}, which allows a much wider class of semigroups to be viewed as ordered categories of various kinds.  (Formal definitions will be given below.)  Further generalisations and specialisations have been studied by numerous authors; of particular relevance to the current work are the restriction semigroups, which themselves trace back to Schweizer and Sklar \cite{SS1960,SS1961,SS1965,SS1967}, and have categorical analogues \cite{CL2002}.  For more details on the significance of such classes of semigroups, and their relation to other parts of algebra and mathematics, see for example \cite{JS2001,Lawson1991,Gould2012,GG2001,BGG2015,KL2017}, and especially~\cite{Gould_notes} for historical background.  For background on semigroups in general, see \cite{Lawson1998,CPbook,Howie,RSbook}.

Recent years have seen a number of important studies of representations of inverse semigroups, especially those of Steinberg  \cite{Steinberg2008,Steinberg2006,Steinberg2010,Steinberg2016}, aspects of which have been extended to Ehresman semigroups by Stein \cite{Stein2016,Stein2019,Stein2017,Stein2017E,Stein2020,MS2020}.  A crucial role in these studies is played by an isomorphism between the semigroup algebra of an appropriate Ehresmann semigroup and an associated category algebra, coming from the Ehresmann-Schein-Nambooripad/Lawson correspondence alluded to above.  Among other consequences, one can then give neat descriptions of simple or indecomposable projective modules.

For the above reasons, it is of considerable interest to find Ehresmann structures on important (families of) semigroups.  Beyond inverse semigroups, classical examples of Ehresmann semigroups include partial transformation semigroups (which are left- but not right-restriction) and monoids of binary relations (which are neither left- nor right-restriction, and not even regular in general).  Roughly speaking, Ehresmann theory is particularly suited to semigroups with some notion of ``partialisation''.  Traditionally, this has meant semigroups whose elements resemble relations or partial transformations in some way (for example, to the authors' knowledge no non-trivial Ehresmann structure is known for the full transformation monoids), but the block symmetry approach of \cite{FL1998} is also fruitful as we will see.

The purpose of the current article is to explore Ehresmann structures on the partition monoids, and related diagram monoids such as (partial) Brauer monoids and rook partition monoids.  An Ehresmann structure is always defined relative to a fixed subsemilattice: i.e., a subsemigroup consisting of commuting idempotents.  Intriguingly, a partition monoid $\P_X$ contains both the symmetric inverse monoid $\I_X$ and the dual symmetric inverse monoid $\J_X$ (all definitions are given below), and therefore contains the semilattices of idempotents coming from both monoids.  One might therefore wonder if both semilattices lead to Ehresmann structures on~$\P_X$; we will see, however, that one semilattice (that of $\J_X$) does lead to such a structure, but the other does not.  Nevertheless, symmetric inverse monoids still do play an extremely important role, leading to a pleasing ``intertwining'' of the above-mentioned partial/block symmetry approach to inverse semigroup theory.  We will also use the general Ehresmann theory to locate some interesting new diagram monoids, which can be thought of as categorical duals to the partial transformation semigroups, in a sense made precise below.  We believe these new monoids are of considerable independent interest.

The paper is organised as follows.  Section \ref{s:pre} contains preliminary background and definitions concerning Ehresmann, inverse and restriction semigroups, various categories and substructures associated to such semigroups, and a discussion of Stein's results on representation theory from~\cite{Stein2017}.  In Section \ref{s:BX} we give an in-depth review of the classical example of monoids of binary relations; the material presented is well-known, even folklore, and is included primarily to set the scene for the following discussion of diagram monoids.  Section \ref{s:PX} concerns the partition monoids~$\P_X$, and contains the main new results of the paper.  We identify two subsemilattices $E$ and $F$ of~$\P_X$---consisting of the idempotents of the submonoids $\I_X$ and $\J_X$, respectively---and show that $\P_X$ is $F$-Ehresmann (Theorem \ref{t:F}) but not $E$-Ehresmann (Remark \ref{r:notE}).  We identify a number of important substructures of $\P_X$ arising from general Ehresmann theory, showing that these are typically isomorphic to symmetric or dual symmetric inverse monoids on generally different base sets (Proposition \ref{p:subs}), or to monoids of various kinds of mappings (e.g., injections between quotients).  We also identify the largest left-, right- and two-sided restriction subsemigroups of $\P_X$ (with respect to $F$).  While the latter is $\J_X$ with an adjoined zero, the other two subsemigroups appear to be new; we will see that they can naturally be thought of as (categorical) duals of the partial transformation semigroups; see Propositions \ref{p:RestPX} and \ref{p:Pfd} and Remark~\ref{r:Pfd}.  Associated categories of (partial) maps between quotients are discussed in Subsection~\ref{ss:CRP}, as well as representations, both of $\P_X$ itself and of the above-mentioned restriction subsemigroups.  Finally, Section \ref{s:ODM} considers other diagram monoids, showing that rook partition monoids have a similar Ehresman structure to ordinary partition monoids, and explaining why the same is not true for Brauer and partial Brauer monoids.  Several ideas for further investigations are discussed throughout the article.

We thank Victoria Gould, Mark Lawson, Stuart Margolis, Itamar Stein and Timothy Stokes for helpful conversations.  We also thank the referee for his/her very careful reading and valuable feedback.

\section{Ehresmann semigroups and restriction semigroups}\label{s:pre}

Here we briefly recall relevant definitions and known results concerning Ehresmann and restriction semigroups.  For the most part we follow \cite{Lawson1991,Stein2017,Gould2012,Gould_notes}, but occasionally adapt notation to suit our current purposes.

\subsection{Ehresmann semigroups}\label{ss:E}

Let $S$ be a semigroup, and write $S^1$ for the monoid obtained by adjoining an identity element if necessary.  Green's $\leqR$, $\leqL$ and $\leqJ$ preorders are defined, for $x,y\in S$, by
\[
x\leqR y \iff x\in yS^1 \COMMA x\leqL y \iff x\in S^1y \AND x\leqJ y \iff x\in S^1yS^1.
\]
So, for example, we have $x\leqR y$ if and only if $x=y$ or else $x=ya$ for some~$a\in S$.  Green's $\R$, $\L$, $\JJ$ and $\H$ relations are then defined by
\[
{\R} = {\leqR}\cap{\geqR} \COMMA 
{\L} = {\leqL}\cap{\geqL} \COMMA 
{\JJ} = {\leqJ}\cap{\geqJ} \AND 
{\H} = {\R}\cap{\L}.
\]
The fifth and final of Green's relations is ${\D}={\R}\vee{\L}$, the join of $\R$ and $\L$ in the lattice of equivalence relations on $S$: i.e., the least equivalence containing both $\R$ and $\L$.  These relations have been of fundamental importance in semigroup theory since their introduction in \cite{Green1951}.  Standard results include that ${\D}={\R}\circ{\L}={\L}\circ{\R}$, and that ${\D}={\JJ}$ if $S$ is finite.  See \cite[Chapter 2]{CPbook} or \cite[Chapter 2]{Howie} for more basic background.

Ehresmann semigroups are defined in terms of certain relations containing $\R$ and $\L$, as we now recall, following the spirit of \cite{Lawson1991,Stein2017}.

Write $E(S)=\set{e\in S}{e=e^2}$ for the set of all idempotents of the semigroup~$S$.  Fix some subset $E\sub E(S)$ that happens to be a semilattice: i.e., $E$ is a commutative subsemigroup of~$S$ consisting entirely of idempotents.  (It is not necessary for $E(S)$ itself to be a semilattice.)  For $x\in S$, we define the (possibly empty) sets
\begin{align}
\nonumber E_L(x) &= \set{e\in E}{x=ex} = \set{e\in E}{x\leqR e} \\
\text{and}\qquad
\label{e:EL_ER} E_R(x) &= \set{e\in E}{x=xe} = \set{e\in E}{x\leqL e},
\end{align}
which consist, respectively, of all left and right identities for $x$ coming from $E$.  The relations $\RE$ and $\LE$ are then defined, for $x,y\in S$, by
\[
x\RE y\iff E_L(x) = E_L(y) \AND x\LE y\iff E_R(x) = E_R(y).
\]
So $x$ and $y$ are $\RE$- or $\LE$-related if and only if they have the same left or right identities from~$E$, respectively.  Clearly $\RE$ and $\LE$ are equivalences.  The relation $\HE$ is defined to be the intersection ${\HE}={\RE}\cap{\LE}$.  It is also easy to check that ${\R}\sub{\RE}$, ${\L}\sub{\LE}$ and ${\H}\sub{\HE}$, where $\R$, $\L$ and $\H$ are the usual Green's relations on $S$, as defined above.

\begin{rem}
%
It is common for the relation $\RE$ to be denoted $\widetilde{\R}^E$ (or similar) in the literature.  It is also common to write $\R^T$ for Green's $\R$ relation on a subsemigroup $T$ of a semigroup $S$ (which is not necessarily the restriction to $T$ of the $\R$ relation on $S$).  However, since we never use the latter notation for subsemigroups (all of Green's relations on a semilattice are trivial), we have opted to write $\RE$ instead of $\widetilde{\R}^E$, omitting the tilde in order to avoid clutter.  Similar comments apply to the $\LE$ and $\HE$ relations.
\end{rem}

With $S$ and $E$ as in the previous paragraph, we say that $S$ is \emph{left $E$-Ehresmann} if the following two conditions hold:
\begin{enumerate}[label=\textup{(L\arabic*)},leftmargin=12mm]
\item \label{L1} Every $\RE$-class contains precisely one element of $E$.
\item \label{L2} $\RE$ is a left congruence on $S$: i.e., $x\RE y\implies ax\RE ay$ for all $a,x,y\in S$.
\een
Dually, $S$ is \emph{right $E$-Ehresmann} if the following two conditions hold:
\begin{enumerate}[label=\textup{(R\arabic*)},leftmargin=12mm]
\item \label{R1} Every $\LE$-class contains precisely one element of $E$.
\item \label{R2} $\LE$ is a right congruence on $S$: i.e., $x\LE y\implies xa\LE ya$ for all $a,x,y\in S$.
\een
The semigroup $S$ is \emph{(two-sided) $E$-Ehresmann} if it is both left and right $E$-Ehresmann.  

Note that all of the above properties are defined with respect to a given semilattice $E\sub E(S)$.
It is also worth noting that any subsemigroup of a left and/or right $E$-Ehresmann semigroup containing~$E$ is still left and/or right $E$-Ehresmann.

Since every element $x$ of an $E$-Ehresmann semigroup satisfies $e\RE x\LE f$ for unique ${e,f\in E}$, it is possible to define two unary operations on $S$ by $x^+=e$ and $x^*=f$.  In fact, $E$-Ehresmann semigroups can be defined equationally as the variety of algebras of type $(2,1,1)$ satisfying various laws, such as $x^+x=x$, $(x^+)^*=x^+$, $(xy)^+=(xy^+)^+$ and more; see \cite[Section~2]{Gould2012}.  For the most part, here we prefer the approach via the relations $\RE$ and $\LE$, as they shed more light on certain natural parameters definable on our monoids of interest.  However, it will be convenient on several occasions to use the ${}^+$ and ${}^*$ operations.

\subsection{Involutions and inverse semigroups}

An \emph{involution} on a semigroup $S$ is a map $S\to S:x\mt x^\circ$ such that $(x^\circ)^\circ=x$ and $(xy)^\circ = y^\circ x^\circ$ for all $x,y\in S$.  It is easy to see that if a left (or right) $E$-Ehresmann semigroup $S$ has an involution ${}^\circ$, then $S$ is also right (or left) $E^\circ$-Ehresmann, respectively, where $E^\circ=\set{e^\circ}{e\in E}$.  Thus, if $E=E^\circ$, then 
\[
[\text{$S$ is left $E$-Ehresmann}] \iff [\text{$S$ is right $E$-Ehresmann}] \iff [\text{$S$ is  $E$-Ehresmann}] .
\]

An important class of involutory semigroups (i.e., semigroups with an involution) are the inverse semigroups.  Recall that a semigroup $S$ is \emph{inverse} \cite{Lawson1998} if for every $x\in S$ there exists a unique element $a\in S$ such that $x=xax$ and $a=axa$; the element $a$ is the \emph{inverse} of $x$, and is generally written as $x^{-1}$.  
It is well known \cite{Schein1963} that a semigroup is inverse if and only if it has a unary operation $x\mt x^\circ({}=x^{-1})$ satisfying the identities
\begin{equation}\label{e:inv}
(x^\circ)^\circ=x=xx^\circ x \COMMA (xy)^\circ=y^\circ x^\circ \AND (xx^\circ)(yy^\circ)=(yy^\circ)(xx^\circ).
\end{equation} 
The idempotents of an inverse semigroup $S$ commute, meaning that $E(S)$ is a semilattice \cite[Theorem~1.17]{CPbook}.  In fact, any inverse semigroup $S$ is $E$-Ehresmann with respect to the entire semilattice $E=E(S)$.  Here ${\RE}={\R}$ and ${\LE}={\L}$ are just the ordinary Green's relations on~$S$, and we have $x^+=xx^{-1}$ and $x^*=x^{-1}x$ for all $x\in S$.

\subsection{Substructures}\label{ss:subs}

If $S$ is $E$-Ehresmann, then we denote the $\HE$-class of $x\in S$ by $H_x^E$, and similarly for $\RE$- and $\LE$-classes.  An important family of submonoids are given by $\HE$-classes of idempotents from~$E$:

\begin{lemma}[{Lawson \cite[Proposition 1.4]{Lawson1991}}]\label{l:HeE}
If $S$ is $E$-Ehresmann, then for any $e\in E$, $H_e^E$ is a monoid with identity $e$.  \epfres
\end{lemma}

Another important subsemigroup of an $E$-Ehresmann semigroup $S$ is defined as follows.  Say $x\in S$ is \emph{$E$-regular} if $x\R e$ and $x\L f$ for some $e,f\in E$.  Since ${\R}\sub{\RE}$ and ${\L}\sub{\LE}$, and since $\RE$- and $\LE$-classes contains unique idempotents from $E$, it follows that $x$ is $E$-regular if and only if $x^+\R x\L x^*$, in the above notation.  We write $\Reg_E(S)=\set{x\in S}{x^+\R x\L x^*}$ for the set of all $E$-regular elemements.  The idea for studying this set goes back at least to \cite{Lawson1990}, in a more general context.

The proof of the following result from \cite{Stein2017} is rather involved, and utilised Prover9 \cite{Prover9}.  

\begin{lemma}[{Stein \cite[Lemma 5.15]{Stein2017}}]\label{l:Reg}
If $S$ is $E$-Ehresmann, then $\Reg_E(S)$ is an inverse subsemigroup of $S$.  \epfres
\end{lemma}

If $S$ is a left $E$-Ehresmann \emph{monoid} with identity $1$, then $E$ must contain $1$ (as the only left or right identity for $1$ is $1$).  In this case, the $\RE$-class $R_1^E$ is a submonoid of $S$.  Dual statements apply to right $E$-Ehresmann monoids.  Further subsemigroups will be discussed in the next subsection.

\subsection{Restriction semigroups}

Recall that a left $E$-Ehresmann semigroup $S$ is \emph{left-restriction} if:
\begin{enumerate}[label=\textup{(L\arabic*)},leftmargin=12mm]\addtocounter{enumi}2
\item \label{L3} For all $x\in S$, $xE\sub Ex$: i.e., for all $x\in S$ and $e\in E$, $xe=fx$ for some $f\in E$.
\een
Dually, a right $E$-Ehresmann semigroup $S$ is \emph{right-restriction} if:
\begin{enumerate}[label=\textup{(R\arabic*)},leftmargin=12mm]\addtocounter{enumi}2
\item \label{R3} For all $x\in S$, $Ex\sub xE$: i.e., for all $x\in S$ and $e\in E$, $ex=xf$ for some $f\in E$.
\een
Using the unary operations ${}^+$ and ${}^*$, these conditions are, respectively, equivalent to:
\begin{enumerate}[label=\textup{(L\arabic*)$'$},leftmargin=12mm]\addtocounter{enumi}2
\item \label{L3'} For all $x\in S$ and $e\in E$, $xe=(xe)^+x$.
\een
\begin{enumerate}[label=\textup{(R\arabic*)$'$},leftmargin=12mm]\addtocounter{enumi}2
\item \label{R3'} For all $x\in S$ and $e\in E$, $ex=x(ex)^*$.
\een
A \emph{restriction semigroup} is one that is both left- and right-restriction.  Again, if $S$ has an involution~${}^\circ$, then $S$ is left-restriction with respect to a semilattice $E$ if and only if it is right-restriction with respect to $E^\circ$.  
And again, any subsemigroup of a left- and/or right-restriction semigroup containing the defining semilattice is still left- and/or right-restriction.

The next result is immediate, and is implicit in \cite[Lemma 3.8]{JS2001}.

\begin{prop}\label{p:Rest}
The largest left-, right- or two-sided restriction subsemigroup of a left-, right- or two-sided $E$-Ehresmann semigroup $S$ is
\begin{align*}
\Rest_L(S,E) &= \set{x\in S}{xE\sub Ex}, \\
\Rest_R(S,E) &= \set{x\in S}{Ex\sub xE}, &\text{or} & &
\Rest(S,E) &= \set{x\in S}{xE= Ex},
\end{align*}
respectively.  \epfres
\end{prop}

\subsection{Ehresmann categories}\label{ss:cats}

Given an $E$-Ehresmann semigroup $S$, Lawson shows how to define a natural category ${C=\bC(S,E)}$ as follows; see \cite[Section 4]{Lawson1991} and \cite[Section 2]{Stein2017}.  The object set of $C$ is $E$.  For $e,f\in E$, the set of morphisms $e\to f$ is
\[
C(e,f) = \set{x\in S}{e\RE x\LE f} = \set{x\in S}{x^+=e,\ x^*=f}.
\]
Composition of morphisms $x\in C(e,f)$ and $y\in C(f,g)$ is simply their product in $S$: i.e., $xy\in C(e,g)$.  (For such $x,y$ we have $y\RE f\implies xy\RE xf = x \RE e$, since $\RE$ is a left congruence, and using $x\LE f$ for $xf=x$; similarly $xy\LE g$, so that indeed $xy\in C(e,g)$.)  
Note that an endomorphism monoid $C(e,e)$, $e\in E$, is precisely the monoid $\HE$-class $H_e^E$.

The category $C = \bC(S,E)$ has the structure of a so-called \emph{Ehresmann category}; see \cite[Section~4]{Lawson1991} or \cite[Definition 2.8]{Stein2017} for the full definition, which we will not reproduce here.  Conversely, an Ehresmann semigroup $\bS(C)$ can be built from an Ehresmann category $C$, and Lawson's main result \cite[Theorem 4.24]{Lawson1991} is that the functors $\bC$ and $\bS$ furnish a category isomorphism between the categories of Ehresmann semigroups and Ehresmann categories.

\subsection{Representations}\label{ss:rep}

The main result of \cite{Stein2017,Stein2017E} concerns another connection between the $E$-Ehresmann semigroup~$S$ and the category~$\bC(S,E)$.  To state it, we also need the partial orders $\leq_r$ and $\leq_l$ defined, respectively, on left and right $E$-Ehresmann semigroups:
\bit
\item If $S$ is left $E$-Ehresmann, then for $x,y\in S$, $x\leq_r y \iff x\in Ey \iff x=x^+y$.
\item If $S$ is right $E$-Ehresmann, then for $x,y\in S$, $x\leq_l y \iff x\in yE \iff x=yx^*$.
\eit
We also say a partial order $\leq$ on a set $X$ is \emph{finite-below} if for every $x\in X$, the set $\set{a\in X}{a\leq x}$ is finite.

Given a small category $C$ and a field $\K$, recall that the \emph{category algebra} $\K[C]$ consists of all $\K$-linear combinations of morphisms from $C$.  The product in $\K[C]$ is first defined for basis elements (morphisms), and then extended by $\K$-linearity.  For morphisms $f$ and $g$, the product~$fg$ is the usual composition in $C$ (another basis element) if defined, or else $fg=0$.

\begin{thm}[{Stein \cite[Theorem 1.5]{Stein2017E}}]\label{t:Stein1}
Let $S$ be an $E$-Ehresmann semigroup that is either left- or right-restriction, and suppose $\leq_r$ or $\leq_l$ (as appropriate) is finite-below.  Then for any field $\K$, the semigroup algebra $\K[S]$ is isomorphic to the category algebra $\K[C]$, where $C=\bC(S,E)$.  \epfres
\end{thm}

The isomorphism $\K[S]\to\K[C]$ maps the basis element $x\in S$ to the (finite) sum $\sum_{a\leq_rx}a$ or $\sum_{a\leq_lx}a$, as appropriate.  The inverse isomorphism involves the M\"obius function of the relevant partial order, but we do not need the details here.  Theorem \ref{t:Stein1} allows the representation theory of the semigroup $S$ to be reduced to that of the category $C$, and in many cases this is a substantial simplification; several examples are given in \cite[Section 5.1]{Stein2017} and also \cite{Stein2016,Stein2019,Steinberg2006,Steinberg2008}.

We also have the following result concerning the maximal semisimple image $\K[S]/\Rad\K[S]$ in the finite case, linking it to the semigroup algebra of the inverse subsemigroup $\Reg_E(S)$ defined in Subsection~\ref{ss:subs}.  An \emph{EI-category} is one in which every endomorphism is an isomorphism.

\begin{thm}[{Stein \cite[Proposition 2.2]{Stein2017E}}]\label{t:Stein2}
Let $S$ be a finite $E$-Ehresmann semigroup that is either left- or right-restriction, and such that $\bC(S,E)$ is an EI-category.  Then for any field $\K$ whose characteristic does not divide the order of any (maximal) subgroup of $S$, $\K[S]/\Rad\K[S]$ is isomorphic to $\K[\Reg_E(S)]$.  \epfres
\end{thm}

The main advantage of Theorem \ref{t:Stein2} is that it allows one to apply the well-developed theory of representations of finite inverse semigroups \cite{Steinberg2006,Steinberg2008,Munn1964}.

\section{Monoids of binary relations and transformations}\label{s:BX}

At this point it is instructive to consider a classical example in depth: namely, monoids of binary relations.  Most, or probably all, of what we say here is well known; indeed, much is considered folklore, and it is not always possible to track down original sources for results.  Many of the ideas are implicit in early work of Tarski on first order logic, as for example in \cite{Tarski1983,Tarski1941,Tarski1955}, though the results on categories and representations come from recent work of Stein \cite{Stein2017,Stein2017E}.  As well as hopefully being a useful and comprehensive source of information on the Ehresmann structure on monoids of binary relations, the material in this section is included to set the scene for the monoids studied in Section \ref{s:PX}, in particular to highlight various parallels and dualities arising.

\subsection{Definitions}\label{ss:defnBX}

Let $X$ be an arbitrary set, and let~$\B_X$ denote the set of all binary relations over~$X$.  Then $\B_X$ is an involutory monoid under composition and converse, defined for $\al,\be\in\B_X$ by
\[
\al\be = \bigset{(x,y)\in X\times X}{(x,u)\in\al,\ (u,y)\in\be\ (\exists u\in X)} \AND \al^\circ=\bigset{(y,x)}{(x,y)\in\al}.
\]
The monoid $\B_X$ is of course well studied in semigroup theory, and in many other parts of mathematics, computer science and logic; see for example \cite{PW1970,Z1963,MP1969,Schein1970,Schein1965,Tarski1941,Tarski1955,Lyndon1950, Lyndon1956,Jonsson1959,HH2001,HJ2012,FS1990,Resende2007,KL2017,BT2016,BT2018,BT2019,BT2020}.

Since several important submonoids will play a key role in all that follows we recall their definitions here as well.  To do so, we first need some terminology.  
The (co)domain and (co)kernel parameters associated to $\al\in\B_X$ are defined by
\begin{align*}
\dom(\al) &= \bigset{x\in X}{(x,u)\in\al\ (\exists u\in X)},\\
\codom(\al) &= \bigset{x\in X}{(u,x)\in\al\ (\exists u\in X)},\\
\ker(\al) &= \bigset{(x,y)\in X\times X}{(x,u),(y,u)\in\al\ (\exists u\in X)},\\
\coker(\al) &= \bigset{(x,y)\in X\times X}{(u,x),(u,y)\in\al\ (\exists u\in X)}.
\end{align*}
So $\dom(\al)$ and $\codom(\al)$ are subsets of $X$, and $\ker(\al)$ and $\coker(\al)$ are relations on $\dom(\al)$ and $\codom(\al)$, respectively; both of these relations are reflexive and symmetric, but need not be transitive.  We say $\al\in\B_X$ is
\bit
\item \emph{injective} if $\ker(\al)$ is the trivial relation on $\dom(\al)$,
\item \emph{coinjective} if $\coker(\al)$ is the trivial relation on $\codom(\al)$,
\item \emph{surjective} if $\codom(\al)=X$,
\item \emph{cosurjective} if $\dom(\al)=X$.
\eit
The submonoids mentioned above are the following:
\bit
\item $\PT_X = \set{\al\in\B_X}{\al\text{ is coinjective}}$, the \emph{partial transformation monoid}, whose elements are partial functions $X\to X$,
\item $\T_X = \set{\al\in\PT_X}{\al\text{ is cosurjective}}$, the \emph{full transformation monoid}, whose elements are (total) functions $X\to X$,
\item $\I_X = \set{\al\in\PT_X}{\al\text{ is injective}}$, the \emph{symmetric inverse monoid}, whose elements are injective partial functions $X\to X$.
\eit
There are also of course anti-isomorphic copies of these: $\PT_X^\circ$, $\T_X^\circ$ and $\I_X^\circ=\I_X$.  Note that~$\I_X$ is inverse (hence its name), with $\al^{-1}=\al^\circ$ for $\al\in\I_X$.  The semilattice $E(\I_X)$ consists of all partial identity maps: $\id_A=\bigset{(a,a)}{a\in A}$, $A\sub X$.  We have $\id_A\id_B=\id_{A\cap B}$ for all $A,B\sub X$.

\subsection{Ehresmann structure}

The idempotents of $\B_X$ itself are rather complicated \cite{MP1969}, but the semilattice
\[
E=E(\I_X)={\set{\id_A}{A\sub X}}
\]
leads to an $E$-Ehresmann structure on $\B_X$, as we now recall; cf.~\cite[Example 4.4]{Stein2017}.  

For $\al\in\B_X$ and $A\sub X$, $\id_A\al=\al\cap(A\times X)$ and $\al\id_A=\al\cap(X\times A)$ are the left- and right-restrictions of $\al$ to $A$, respectively.  It quickly follows that
\[
E_L(\al) = \set{\id_A}{A\supseteq\dom(\al)} \AND E_R(\al) = \set{\id_A}{A\supseteq\codom(\al)} \qquad\text{for all $\al\in\B_X$.}
\]
Thus,
\[
\al\RE\be \iff \dom(\al)=\dom(\be) \AND \al\LE\be \iff \codom(\al)=\codom(\be).
\]
An $\RE$-class is determined by the common domain of its elements, so contains a unique element of $E$; thus \ref{L1} holds.  If $\al,\th\in\B_X$, then $\dom(\th\al) = \bigset{x\in X}{(x,u)\in\th\ (\exists u\in\dom(\al))}$.  It follows that $\dom(\al)=\dom(\be)\implies\dom(\th\al)=\dom(\th\be)$ for all $\al,\be,\th\in\B_X$; thus, \ref{L2} holds.  Since $E=E^\circ$, \ref{R1} and \ref{R2} hold as well.  

This all shows that $\B_X$ is $E$-Ehresmann.  It is worth noting that for $\al\in\B_X$ we have ${\al^+=\id_{\dom(\al)}}$ and ${\al^*=\id_{\codom(\al)}}$.  It quickly follows that the $\leq_r$ relation has a simple characterisation.  Namely, given $\al,\be\in\B_X$ we have $\al\leq_r\be$ if and only if $\al$ is the restriction of $\be$ to~$\dom(\al)$.  A dual statement holds for the $\leq_l$ relation.

\subsection{Substructures}

The monoid $\HE$-class of a partial identity $\id_A\in E$ (cf.~Lemma \ref{l:HeE}) is
\[
H_{\id_A}^E = \set{\al\in\B_X}{\dom(\al)=\codom(\al)=A},
\]
the monoid of all surjective and cosurjective binary relations on $A$.  

Although Green's relations on $\B_X$ are again fairly complicated in general \cite{PW1970,Z1963}, the inverse subsemigroup
\[
\Reg_E(\B_X) = \set{\al\in\B_X}{\al^+\R\al\L\al^*}
\]
consisting of all $E$-regular elements (cf.~Lemma \ref{l:Reg}) can be readily described as follows.  Consider some $\al\in\Reg_E(\B_X)$, and write $A=\dom(\al)$ and $B=\codom(\al)$.  Keeping in mind that $\al^+=\id_A$ and $\al^*=\id_B$, we have
\[
\id_A = \al\rho \AND \id_B = \lam\al \qquad\text{for some $\lam,\rho\in\B_X$.}
\]
Since $\al\rho=(\al\id_B)\rho=\al(\id_B\rho)$, we may replace $\rho$ by $\id_B\rho$ to assume that $\dom(\rho)\sub B$, and similarly that $\codom(\lam)\sub A$.  But then $\lam = \lam\id_A = \lam(\al\rho) = (\lam\al)\rho = \id_B\rho = \rho$.  Since then $\id_A = \al\rho$ and $\id_B = \rho\al$, it follows that $\al$ is a bijection $A\to B$ (with inverse $\rho$).  This all shows that
\[
\Reg_E(\B_X) = \I_X
\]
is precisely the symmetric inverse monoid.

Since $\B_X$ is a monoid, we may also consider the $\RE$- and $\LE$-classes of the identity $\id_X$.  These are the submonoids
\begin{align*}
R_{\id_X}^E &= \set{\al\in\B_X}{\dom(\al)=X} &&\text{and}& L_{\id_X}^E &= \set{\al\in\B_X}{\codom(\al)=X}\\
 &= \set{\al\in\B_X}{\al\text{ is cosurjective}} &&&  &= \set{\al\in\B_X}{\al\text{ is surjective}}.
\end{align*}

\subsection{Restriction subsemigroups}

The monoid $\B_X$ is neither left- nor right-restriction for $|X|\geq2$, as follows from the easily checked fact that
\[
\text{$\al E\sub E\al \iff \al$ is coinjective} \AND \text{$E\al\sub \al E \iff \al$ is injective} \qquad\text{for all $\al\in\B_X$.}
\]
It follows immediately (cf.~Proposition \ref{p:Rest}) that the largest left-restriction subsemigroup of~$\B_X$ is
\[
\Rest_L(\B_X,E) = \set{\al\in\B_X}{\al\text{ is coinjective}}  = \PT_X ,
\]
the partial transformation monoid.  Similarly, $\Rest_R(\B_X,E) = \PT_X^\circ$ is an anti-isomorphic copy of $\PT_X$.  The largest restriction subsemigroup of $\B_X$ is
\[
\Rest(\B_X,E) = \PT_X\cap\PT_X^\circ =  \set{\al\in\B_X}{\al\text{ is injective and coinjective}} = \I_X,
\]
the symmetric inverse monoid.

The partial transformation monoid $\PT_X$ is of course an $E$-Ehresmann semigroup in its own right, with respect to the same semilattice $E=E(\I_X)$.  In the semigroup $\PT_X$, the monoid $\HE$-class of $\id_A$, $A\sub X$, is the set
\[
H_{\id_A}^E = \set{\al\in\PT_X}{\dom(\al)=\codom(\al)=A},
\]
which this time is the monoid of all surjective (total) functions $A\to A$; when $A$ is finite, this is of course the \emph{symmetric group} $\S_A$, consisting of all permutations of $A$, and coinciding in this (finite) case with the $\H$-class of $\id_A$ in $\PT_X$.  The inverse subsemigroup $\Reg_E(\PT_X)$ is still the symmetric inverse monoid $\I_X$.  The $\RE$- and $\LE$-classes in $\PT_X$ of the identity $\id_X$ are
\[
R_{\id_X}^E = \set{\al\in\PT_X}{\dom(\al)=X} \AND L_{\id_X}^E = \set{\al\in\PT_X}{\codom(\al)=X}.
\]
The first of these is precisely the full transformation monoid $\T_X$, while the second is the monoid of all surjective partial functions $X\to X$.  It is easily seen that the latter is precisely the $\L$-class in $\B_X$ of the identity $\id_X$: i.e., the set of all left-invertible relations from $\B_X$.  When $X$ is finite, $L_{\id_X}^E$ is the symmetric group $\S_X$ (but $R_{\id_X}^E$ is not, as we have already observed that $R_{\id_X}^E=\T_X$ for any $X$).

\subsection{Categories and representations}\label{ss:catBX}

Identifying an idempotent $\id_A\in E$ with the subset $A\sub X$, the category $C = \bC(\B_X,E)$ has object set $\Sub(X) = \set{A}{A\sub X}$.  The morphism set
\[
C(A,B) = \set{\al\in\B_X}{\id_A \RE \al \LE \id_B} = \set{\al\in\B_X}{\dom(\al)=A,\ \codom(\al)=B}
\]
consists of all surjective and cosurjective relations $A\to B$.  For more on representations on monoids of binary relations see \cite{BT2016,BT2018,BT2019,BT2020}.

The category $C'=\bC(\PT_X,E)$ has the same object set, but this time the morphism set
\[
C'(A,B) = \set{\al\in\PT_X}{\dom(\al)=A,\ \codom(\al)=B}
\]
consists of all surjective (total) functions $A\to B$.  When $X$ is finite, the partial order $\leq_r$ is certainly finite-below, so Theorem \ref{t:Stein1} applies, and reduces the representation theory of $\PT_X$ to that of the category $C'=\bC(\PT_X,E)$.  This has been explored in depth by Stein in \cite{Stein2016,Stein2019}.  

For finite $A\sub X$, the endomorphism monoid $C'(A,A)$ is precisely the symmetric group $\S_A$.  It follows that $C'=\bC(\PT_X,E)$ is an EI-category when $X$ is finite.  Thus, Theorem~\ref{t:Stein2} applies, and tells us that for any field $\K$ of suitable characteristic (either $0$ or else greater than~$|X|$, noting that maximal subgroups of $\PT_X$ are symmetric groups of degrees $0,1,\ldots,n$), the maximal semisimple image of $\K[\PT_X]$ is precisely~$\K[\I_X]$.

\section{Partition monoids}\label{s:PX}

We now turn our attention to the partition monoid $\P_X$.  After recalling the definitions, we observe that $\P_X$ contains two natural semilattices, and show that only one of these leads to an Ehresmann structure.  We then locate the various substructures and categories arising from the general theory, identify the largest left- and/or right-restriction subsemigroups, and make some comments on representation theory.  

Along the way, we uncover some striking parallels and dualities between $\B_X$ (cf.~Section~\ref{s:BX}) and~$\P_X$.  For example, full transformation monoids and symmetric inverse monoids arise naturally, but so too does the \emph{dual} symmetric inverse monoid.  The role played in Section~\ref{s:BX} by surjections between subsets of $X$ is here played by injections between quotients of~$X$.  Other parallels between $\B_X$ and $\P_X$ have been explored in \cite{FL2011,FitzGerald2013}.

\subsection{Definitions}\label{ss:definitions}

Fix a set $X$, and let $X'=\set{x'}{x\in X}$ be a disjoint copy of $X$.  The partition monoid $\P_X$ consists of all set partitions of $X\cup X'$, under a product defined shortly.  When $X=\{1,\ldots,n\}$ for some positive integer $n$, we write $\P_n$ for $\P_X$.  A partition from $\P_X$ will be identified with any graph on vertex set $X\cup X'$ whose connected components are the blocks of the partition; when depicting such a graph, vertices from $X$ and $X'$ are displayed on upper and lower rows, respectively.  For example, the partitions
\begin{align}
\nonumber \al &= \big\{ \{1,4\},\{2,3,4',5'\},\{5,6\},\{1',2',6'\},\{3'\}\big\} \\
\label{e:albe} \text{and}\qquad \be &= \big\{ \{1,2\}, \{3,4,1'\}, \{5,4',5',6'\}, \{6\}, \{2'\}, \{3'\} \big\}
\end{align}
from $\P_6$ are shown in Figure \ref{f:P6}.

Partitions $\al,\be\in\P_X$ are multiplied as follows.  First, let $\al_\downarrow$ and $\be^\uparrow$ be the graphs obtained by changing every lower vertex $x'$ of $\al$, and every upper vertex $x$ of $\be$, to $x''$.  Now let $\Pi(\al,\be)$ be the graph on vertex set $X\cup X''\cup X'$ with edge set the union of the edge sets of $\al_\downarrow$ and $\be^\uparrow$;~$X''$ is the middle row of $\Pi(\al,\be)$.  The product $\al\be\in\P_X$ is then the partition of $X\cup X'$ such that $u,v\in X\cup X'$ belong to the same block of $\al\be$ if and only if there is a path from $u$ to $v$ in~$\Pi(\al,\be)$.  For example, with $\al,\be\in\P_6$ as in \eqref{e:albe}, Figure \ref{f:P6} shows how to calculate the product 
\[
\al\be = \big\{\{1,4\}, \{2,3,1',4',5',6'\}, \{5,6\} , \{2'\}, \{3'\} \big\}.
\]
The operation is associative, so $\P_X$ is a monoid with identity $\id_X={\bigset{\{x,x'\}}{x\in X}}$.

\begin{figure}[h]
\begin{center}
\begin{tikzpicture}[scale=.5]

\begin{scope}[shift={(0,0)}]	
\uvs{1,...,6}
\lvs{1,...,6}
\uarcx14{.6}
\uarcx23{.3}
\uarcx56{.3}
\darc12
\darcx26{.6}
\darcx45{.3}
\stline34
\draw(0.6,1)node[left]{$\alpha=$};
\draw[->](7.5,-1)--(9.5,-1);
\end{scope}

\begin{scope}[shift={(0,-4)}]	
\uvs{1,...,6}
\lvs{1,...,6}
\uarc12
\uarc34
\darc45
\darc56
\stline31
\stline55
\draw(0.6,1)node[left]{$\beta=$};
\end{scope}

\begin{scope}[shift={(10,-1)}]	
\uvs{1,...,6}
\lvs{1,...,6}
\uarcx14{.6}
\uarcx23{.3}
\uarcx56{.3}
\darc12
\darcx26{.6}
\darcx45{.3}
\stline34
\draw[->](7.5,0)--(9.5,0);
\end{scope}

\begin{scope}[shift={(10,-3)}]	
\uvs{1,...,6}
\lvs{1,...,6}
\uarc12
\uarc34
\darc45
\darc56
\stline31
\stline55
\end{scope}

\begin{scope}[shift={(20,-2)}]	
\uvs{1,...,6}
\lvs{1,...,6}
\uarcx14{.6}
\uarcx23{.3}
\uarcx56{.3}
\darc14
\darc45
\darc56
\stline21
\draw(6.4,1)node[right]{$=\alpha\beta$};
\end{scope}

\end{tikzpicture}
\caption{Multiplication of the two partitions $\al,\be\in\P_6$ given in \eqref{e:albe}.}
\label{f:P6}
\end{center}
\end{figure}

At times it will be convenient to use a standard tabular notation for partitions.  First, we say that a block of a partition $\al\in\P_X$ is a \emph{transversal} if it contains points from both $X$ and~$X'$.  A block of $\al$ is an \emph{upper} or \emph{lower non-transversal} if it is contained in $X$ or $X'$, respectively.  For $\al\in\P_X$ we write
\[
\al = \begin{partn}{2} A_i&C_j\\ \hhline{~|-} B_i&D_k\end{partn}_{i\in I,\ j\in J,\ k\in K}
\]
to indicate that $\al$ has transversals $A_i\cup B_i'$ ($i\in I$), upper non-transversals $C_j$ ($j\in J$), and lower non-transversals $D_k'$ ($k\in K$).  This will often be abbreviated to $\al = \begin{partn}{2} A_i&C_j\\ \hhline{~|-} B_i&D_k\end{partn}$, with the indexing sets $I$, $J$ and $K$ being implied, rather than explicitly named.  For example, with $\al,\be\in\P_6$ as in~\eqref{e:albe}, we have
\[
\al = \begin{partn}{3} 2,3&1,4&5,6\\ \hhline{~|-|-} 4,5&1,2,6&3\end{partn}
\AND
\be = \begin{partn}{4} 3,4&5&1,2&6\\ \hhline{~|~|-|-} 1&4,5,6&2&3\end{partn}.
\]

The monoid $\P_X$ has an involution $\al\mt\al^\circ$, where $\al^\circ$ is obtained from $\al$ by interchanging dashed and undashed elements of $X\cup X'$ (equivalently, by reflecting graphs in a horizontal axis).  In tabular notation, $\begin{partn}{2} A_i&C_j\\ \hhline{~|-} B_i&D_k\end{partn}^\circ=\begin{partn}{2} B_i&D_k\\ \hhline{~|-} A_i&C_j\end{partn}$.
In fact, we have
\[
(\al^\circ)^\circ = \al = \al\al^\circ\al \AND (\al\be)^\circ=\be^\circ\al^\circ \qquad\text{for all $\al,\be\in\P_X$.}
\]
These identities show that $\P_X$ is a \emph{regular $*$-monoid} in the sense of Nordahl and Scheiblich \cite{NS1978}.  Since the identity $(\al\al^\circ)(\be\be^\circ)=(\be\be^\circ)(\al\al^\circ)$ does not hold for $|X|\geq2$, it follows that $\P_X$ is not inverse; cf.~\eqref{e:inv}.

\subsection{Green's relations and (inverse) submonoids}\label{ss:Green}

Green's relations on $\P_X$ may be characterised in terms of the (co)domain, (co)kernel and rank parameters, defined for $\al\in\P_X$ as follows:
\begin{align*}
\dom(\al) &= \set{x\in X}{x\text{ belongs to a transversal of }\al},\\
\codom(\al) &= \set{x\in X}{x'\text{ belongs to a transversal of }\al},\\
\ker(\al) &= \set{(x,y)\in X\times X}{\text{$x$ and $y$ belong to the same block of }\al},\\
\coker(\al) &= \set{(x,y)\in X\times X}{\text{$x'$ and $y'$ belong to the same block of }\al},\\
\rank(\al) &= \text{the number of transversals of $\al$}.
\end{align*}
For example, $\al$ in \eqref{e:albe} has rank equal to $1$, domain equal to $\{2,3\}$, and its cokernel-classes are $\{1,2,6\}$, $\{3\}$ and $\{4,5\}$.  As with binary relations (cf.~Subsection \ref{ss:defnBX}), $\dom(\al)$ and $\codom(\al)$ are subsets of $X$; however, note that $\ker(\al)$ and $\coker(\al)$ are always equivalences on $X$, not just relations on $\dom(\al)$ and $\codom(\al)$.  Also, $\rank(\al)$ is a cardinal between $0$ and $|X|$.

The next result follows from \cite[Lemma 3.1]{FL2011}:

\begin{prop}\label{p:RL}
For $\al,\be\in\P_X$, we have
\ben
\item \label{i:RL1} $\al\leqR\be \iff [\ker(\be)\sub\ker(\al)$ and every upper non-transversal of $\be$ is a block of $\al]$,
\item \label{i:RL2} $\al\leqL\be \iff [\coker(\be)\sub\coker(\al)$ and every lower non-transversal of $\be$ is a block of~$\al]$,
\item \label{i:RL5} $\al\leqJ\be \iff \rank(\al)\leq\rank(\be)$,
\item \label{i:RL3} $\al\R\be \iff [\dom(\al)=\dom(\be)$ and $\ker(\al)=\ker(\be)]$,
\item \label{i:RL4} $\al\L\be \iff [\codom(\al)=\codom(\be)$ and $\coker(\al)=\coker(\be)]$,
\item \label{i:RL6} $\al\JJ\be \iff \al\D\be \iff \rank(\al)=\rank(\be)$.  \epfres
\een
\end{prop}

The above parameters also allow us to conveniently describe certain natural submonoids that will be important in all that follows.  Writing $\De=\bigset{(x,x)}{x\in X}$ for the trivial equivalence, we define
\bit
\item $\I_X = \set{\al\in\P_X}{\ker(\al)=\coker(\al)=\De}$,
\item $\J_X = \set{\al\in\P_X}{\dom(\al)=\codom(\al)=X}$,
\item $\T_X = \set{\al\in\P_X}{\dom(\al)=X,\ \coker(\al)=\De}$.
\eit
As noted in \cite[Section 2]{EF2012}, these submonoids are isomorphic to the \emph{symmetric inverse monoid}, the \emph{dual symmetric inverse monoid} and the \emph{full transformation monoid} over $X$, justifying our (re)use of the notation $\I_X$ and $\T_X$.  (There is also an anti-isomorphic copy $\T_X^\circ$ of $\T_X$.)  As in Subsection~\ref{ss:defnBX}, the elements of $\I_X$ and $\T_X$ are (identified with) partial bijections $X\to X$ or (total) functions~${X\to X}$.  An element of $\J_X$ is (identified with) a \emph{block bijection} on $X$: i.e., a bijection $\bA\to\bB$ for partitions $\bA,\bB$ of $X$ \cite{FL1998}.  See Figure \ref{f:I6J6} for some examples; the pictured block bijection maps~${\{1,2\}\mt\{1,2,3\}}$, $\{3,5\}\mt\{4\}$ and $\{4,6\}\mt\{5,6\}$.

\begin{figure}[h]
\begin{center}
\begin{tikzpicture}[scale=.5]

\begin{scope}[shift={(0,0)}]	
\uvs{1,...,6}
\lvs{1,...,6}
\stline13
\stline21
\stline44
\stline56
\end{scope}

\begin{scope}[shift={(10,0)}]	
\uvs{1,...,6}
\lvs{1,...,6}
\stline11
\stline34
\stline54
\stline45
\stline66
\uarc12
\darc12
\darc23
\darc56
\end{scope}

\begin{scope}[shift={(20,0)}]	
\uvs{1,...,6}
\lvs{1,...,6}
\stline13
\stline21
\stline34
\stline44
\stline56
\stline64
\end{scope}

\end{tikzpicture}
\caption{Sample elements of $\I_6$ (left), $\J_6$ (middle) and $\T_6$ (right).}
\label{f:I6J6}
\end{center}
\end{figure}
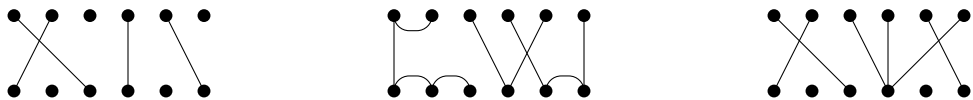

\begin{rem}\label{r:part_bb}
In fact, any element $\al$ of $\P_X$ can be thought of as a \emph{partial block bijection} of~$X$: i.e., a partial bijection between the quotients $X/\ker(\al)\to X/\coker(\al)$.  (For an equivalence $\ve$ on $X$ we write $X/\ve$ for the associated partition of $X$, consisting of all $\ve$-classes.)  The mapping rule is determined by the transversals: $A\in X/\ker(\al)$ maps to $B\in X/\coker(\al)$ if and only if $A\cup B'$ is a transversal of $\al$.  So, for example, $\al\in\P_X$ from Figure \ref{f:P6} is identified with the partial bijection
\[
\big\{ \{1,4\}, \{2,3\}, \{5,6\}\big\} \to \big\{ \{1,2,6\}, \{3\}, \{4,5\} \big\}
\]
under which $\{2,3\}$ is mapped to $\{4,5\}$.  This perspective will be a recurring theme in what follows.
\end{rem}

Although $\P_X$ is not itself inverse, the submonoids~$\I_X$ and $\J_X$ both are; in both submonoids, the unique inverse of $\al$ is $\al^\circ$, defined at the end of Subsection \ref{ss:definitions}.  In particular, the sets $E(\I_X)$ and $E(\J_X)$ are both semilattices.  Since these are both contained in $\P_X$, they are both candidates for Ehresmann structures.  We explore this in the next subsection.

We write $\Sub(X)$ and $\Eq(X)$ for the sets of all subsets of $X$ and all equivalences on $X$, respectively.  For $A\in\Sub(X)$ and $\ve\in\Eq(X)$, we define the partitions
\[
\id_A = \bigset{\{a,a'\}}{a\in A} \cup \bigset{\{x\},\{x'\}}{x\in X\sm A} \AND
\id_\ve = \set{A\cup A'}{A\in X/\ve}.
\]
Examples are given in Figure \ref{f:EF}.  As in \cite{Lipscombe1996,FL1998}, we have
\[
E(\I_X) = \bigset{\id_A}{A\in\Sub(X)} \AND E(\J_X) = \bigset{\id_\ve}{\ve\in\Eq(X)}.
\]
We noted earlier that $\id_A\id_B=\id_{A\cap B}$ for ${A,B\in\Sub(X)}$.  On the other hand, ${\id_\ve\id_\eta=\id_{\ve\vee\eta}}$ for $\ve,\eta\in\Eq(X)$, where $\ve\vee\eta$ is the \emph{join} of~$\ve$ and~$\eta$: i.e., the least equivalence containing $\ve\cup\eta$.

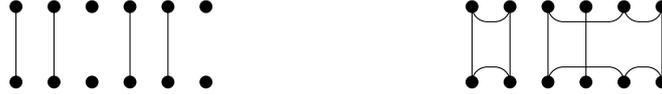
\begin{figure}[h]
\begin{center}
\begin{tikzpicture}[scale=.5]

\begin{scope}[shift={(0,0)}]	
\uvs{1,...,6}
\lvs{1,...,6}
\stline11
\stline22
\stline44
\stline55
\end{scope}

\begin{scope}[shift={(12,0)}]	
\uvs{1,...,6}
\lvs{1,...,6}
\stline11
\stline22
\stline33
\stline44
\stline66
\uarc12
\darc12
\uarc35
\uarc56
\darc35
\darc56
\end{scope}

\end{tikzpicture}
\caption{Idempotents $\id_A\in \I_6$ (left) and $\id_\ve\in \J_6$ (right), where $A=\{1,2,4,5\}$ and $\ve$ has classes $\{1,2\}$, $\{3,5,6\}$ and $\{4\}$.}
\label{f:EF}
\end{center}
\end{figure}

\subsection{Ehresmann structure}\label{ss:Ehresmann}

As noted in the previous subsection, we wish to use the semilattices of idempotents from the inverse monoids $\I_X$ and $\J_X$ to investigate possible Ehresmann structures on the partition monoid~$\P_X$.  Accordingly, we define
\[
E = E(\I_X) \AND F = E(\J_X).
\]
The first step is to identify the sets $E_L(\al)$, $F_L(\al)$, and so on.  This requires the (co)support parameters associated to a partition $\al\in\P_X$:
\begin{align*}
\supp(\al) &= \set{x\in X}{\{x\}\text{ is not a block of $\al$}}\\
\cosupp(\al) &= \set{x\in X}{\{x'\}\text{ is not a block of $\al$}}.
\end{align*}
For example, $\be\in\P_6$ from Figure \ref{f:P6} has support $\{1,2,3,4,5\}$ and cosupport $\{1,4,5,6\}$.  Note that $\supp(\al)=\cosupp(\al)=X$ for all $\al\in\J_X$, while $\supp(\be)=\dom(\be)$ and $\cosupp(\be)=\codom(\be)$ for all $\be\in\I_X$.  

\begin{lemma}\label{l:eq}
For $\al\in\P_X$ we have
\ben\bmc2
\item \label{i:eq1} $E_L(\al) = \bigset{\id_A}{A\supseteq\supp(\al)}$,
\item \label{i:eq2} $E_R(\al) = \bigset{\id_A}{A\supseteq\cosupp(\al)}$,
\item \label{i:eq3} $F_L(\al) = \bigset{\id_\ve}{\ve\sub\ker(\al)}$,
\item \label{i:eq4} $F_R(\al) = \bigset{\id_\ve}{\ve\sub\coker(\al)}$.
\emc\een
\end{lemma}

\pf
We just prove \ref{i:eq1} and \ref{i:eq3}, as the other statements are dual.

\pfitem{\ref{i:eq1}}  Let $A\in\Sub(X)$.  By definition (cf.~\eqref{e:EL_ER}) we have $\id_A\in E_L(\al) \iff \al\leqR\id_A$.  Since $\id_A$ has trivial kernel, the latter is equivalent (using Proposition \ref{p:RL}\ref{i:RL1}) to every upper non-transversal of $\id_A$ being a block of $\al$, which is equivalent to $\supp(\al)\sub A$.  

\pfitem{\ref{i:eq3}}  Let $\ve\in\Eq(X)$.  By definition (cf.~\eqref{e:EL_ER}) we have $\id_\ve\in F_L(\al) \iff \al\leqR\id_\ve$, so the statement now follows from Proposition \ref{p:RL}\ref{i:RL1}, given that $\id_\ve$ has no non-transversals.
\epf

As an immediate consequence, the relations $\RE$, $\LE$, $\RF$ and $\LF$ have the following descriptions:

\begin{prop}\label{p:RELEPX}
For $\al,\be\in\P_X$, and with $E=E(\J_X)$ and $F=E(\I_X)$ as above, we have
\ben\bmc2
\item $\al \RE \be \iff \supp(\al)=\supp(\be)$,
\item $\al \LE \be \iff \cosupp(\al)=\cosupp(\be)$,
\item $\al \RF \be \iff \ker(\al)=\ker(\be)$,
\item $\al \LF \be \iff \coker(\al)=\coker(\be)$.  \epfres
\emc\een
\end{prop}

It quickly follows that \ref{R1} and \ref{L1} hold for both semilattices $E$ and $F$.  For example, the unique elements of $E$ and $F$ in the $\RE$- and $\RF$-classes of $\al$ are $\id_{\supp(\al)}$ and $\id_{\ker(\al)}$.

\begin{lemma}
The equivalence $\RF$ is a left congruence on $\P_X$, and $\LF$ is a right congruence.
\end{lemma}

\pf
We just prove the statement concerning $\RF$, as the other is dual.  
To do so, suppose $\al,\be,\th\in\P_X$ and $\al\RF\be$: i.e., $\ker(\al)=\ker(\be)$.  We must show that $\ker(\th\al)=\ker(\th\be)$, and by symmetry it suffices to show that $\ker(\th\al)\sub\ker(\th\be)$.  

Suppose $(x,y)\in\ker(\th\al)$.  Then there is a path $\pi$ from $x$ to $y$ in the graph $\Pi(\th,\al)$, as defined in Subsection~\ref{ss:definitions}.  If $\pi$ stays entirely in the ``$\th$-half'' of $\Pi(\th,\al)$, then this path is present in the graph $\Pi(\th,\be)$ as well, so that $(x,y)\in\ker(\th\be)$ in this case.  Otherwise, $\pi$ has the form
\begin{equation}\label{e:path}
x \lar\th z_1'' \lar\al z_2'' \lar\th z_3'' \lar\al \cdots \lar\th z_{k-1}'' \lar\al z_k'' \lar\th y \qquad\text{for some $z_1,\ldots,z_k\in X$.}
\end{equation}
Here, for example, $z_1''\lar\al z_2''$ indicates a path from $z_1''$ to $z_2''$ contained in the ``$\al$-half'' of $\Pi(\th,\al)$.  Since then $(z_1,z_2),(z_3,z_4),\ldots,(z_{k-1},z_k)\in\ker(\al)=\ker(\be)$, it follows that there is a path of the form \eqref{e:path} in the graph $\Pi(\th,\be)$, so again we have $(x,y)\in\ker(\th\be)$.  
\epf

It follows that conditions \ref{R2} and \ref{L2} hold for the semilattice $F$, so we have the following:

\begin{thm}\label{t:F}
For any set $X$, the partition monoid $\P_X$ is $F$-Ehresmann with respect to the semilattice ${F=E(\J_X)}$, where $\J_X$ is the dual symmetric inverse monoid.  \epfres
\end{thm}

\begin{rem}\label{r:notE}
Neither \ref{R2} nor \ref{L2} hold for the semilattice $E=E(\I_X)$ if $|X|\geq2$.  For example, consider the partitions $\al = \custpartn{1,2}{1,2}{\uarc12 \darc12}$, $\be = \custpartn{1,2}{1,2}{\stline11\stline22}$ and $\th=\custpartn{1,2}{1,2}{\stline11}$, all from $\P_2$.  Then $\supp(\al)=\supp(\be)$ and $\cosupp(\al)=\cosupp(\be)$, yet $\supp(\th\al)\not=\supp(\th\be)$ and $\cosupp(\al\th)\not=\cosupp(\be\th)$.  Thus,~$\P_X$ is not $E$-Ehresmann.
\end{rem}

We conclude this subsection with a description of the $\leq_r$ and $\leq_l$ partial orders associated to the $F$-Ehresmann structure on $\P_X$.  Recall from Subsection \ref{ss:rep} that these are defined, for $\al,\be\in\P_X$, by
\[
\al\leq_r\be \iff \al\in F\be \AND \al\leq_l\be \iff \al\in \be F.
\]
For the following statement, if $\al,\be\in\P_X$, we write $\al\pre\be$ to indicate that $\al$ \emph{refines} $\be$: i.e., that every block of $\al$ is contained in a block of $\be$, or equivalently if every block of $\be$ is a union of blocks of $\al$.

\begin{prop}\label{p:leqr}
For $\al,\be\in\P_X$, we have
\ben
\item \label{i:leqr1} $\al\leq_r\be \iff [\be\pre\al$ and every lower non-transversal of $\be$ is a block of $\al]$,
\item \label{i:leqr2} $\al\leq_l\be \iff [\be\pre\al$ and every upper non-transversal of $\be$ is a block of $\al]$.
\een
\end{prop}

\pf
As usual, it suffices to prove \ref{i:leqr1}.  Beginning with the forwards implication, suppose ${\al\leq_r\be}$, so that $\al=\id_\ve\be$ for some $\ve\in\Eq(X)$.  Write $\be = \begin{partn}{2} A_i&C_j\\ \hhline{~|-} B_i&D_k\end{partn}$.  Every lower non-transversal~$D_k'$ survives in the product $\id_\ve\be=\al$, so it remains to show that $\be\pre\al$.  Next note that
\[
\al = \id_\ve\be = \id_\ve(\id_{\ker(\be)}\be) = \id_{\ve\vee\ker(\be)}\be,
\]
so we may assume without loss of generality that every $\ve$-class is a union of $\ker(\be)$-classes.  Write $X/\ve=\set{G_l}{l\in L}$.  The condition on $\ve$- and $\ker(\be)$-classes tells us that for each $l\in L$ there exist subsets $I_l\sub I$ and $J_l\sub J$ such that $G_l = \bigcup_{i\in I_l}A_i \cup \bigcup_{j\in J_l}C_j$.  Note that one of $I_l$ or $J_l$ may be empty, but not both.  Note also that $I=\bigcup_{l\in L}I_l$ and $J=\bigcup_{l\in L}J_l$, and that the $I_l$ are pairwise-disjoint, as are the $J_l$.  But then for each $l$, $\al=\id_\ve\be$ contains the block $\bigcup_{i\in I_l}(A_i\cup B_i') \cup \bigcup_{j\in J_l}C_j$, which is a union of blocks of $\al$ (this block is a transversal if and only if $I_l$ is non-empty).  Together with the $D_k$, these are all the blocks of $\id_\ve\be=\al$, so it follows that $\be\pre\al$.

For the converse, suppose now that $\be\pre\al$ and that every lower non-transversal of $\be$ is a block of $\al$.  Again write $\be = \begin{partn}{2} A_i&C_j\\ \hhline{~|-} B_i&D_k\end{partn}$.  So each $D_k'$ is a block of $\al$; let the remaining blocks of $\al$ be $\set{H_l}{l\in L}$.  Since $\be\pre\al$, each $H_l$ has the form $H_l = \bigcup_{i\in I_l}(A_i\cup B_i') \cup \bigcup_{j\in J_l}C_j$ for some $I_l\sub I$ and $J_l\sub J$.  It then quickly follows that $\al=\id_\ve\be$, where the $\ve$-classes are $\bigcup_{i\in I_l}A_i \cup \bigcup_{j\in J_l}C_j$ for each $l\in L$.
\epf

\begin{rem}\label{r:leqr}
The relations $\leq_r$ and $\leq_l$ have natural geometrical interpretations.  For example, given a partition $\be=\begin{partn}{2} A_i&C_j\\ \hhline{~|-} B_i&D_k\end{partn}$ from $\P_X$, we can obtain any $\al\leq_r\be$ by keeping the lower non-transversals $D_k'$, and ``fusing'' together any of the remaining blocks of $\be$ into larger blocks.  Figure~\ref{f:leqr} shows this process, illustrating a product $\al=\id_\ve\be$, where (as in the proof of Proposition \ref{p:leqr}) every $\ve$-class is a union of $\ker(\be)$-classes.
\end{rem}

\begin{figure}[h]
\begin{center}
\begin{tikzpicture}[scale=.85]

\begin{scope}[shift={(0,2)}]	
\colrect1717{blue}
\colrect{7.5}{10.5}{7.5}{10.5}{blue}
\colrect{11}{18}{11}{18}{blue}
\draw[|-] (.5,2)--(.5,0);
\node[left] () at (.5,1) {$\id_\ve$};
\end{scope}

\begin{scope}[shift={(0,0)}]	
\colrect{3.5}7{2.5}5{red}
\colrect{11}{13}{12}{13}{red}
\colrect{13.5}{14.5}{13.5}{15}{red}
\colrect{15}{18}{15.5}{16}{red}
\colrecthigh1{1.5}{red}
\colrecthigh23{red}
\colrecthigh{7.5}{8.5}{red}
\colrecthigh9{9.5}{red}
\colrecthigh{10}{10.5}{red}
\colrectlow9{10}{red}
\colrectlow{5.5}{6.5}{red}
\colrectlow12{red}
\colrectlow7{7.5}{red}
\colrectlow8{8.5}{red}
\colrectlow{10.5}{11.5}{red}
\colrectlow{16.5}{18}{red}
\draw[|-|] (.5,2)--(.5,0);
\node[left] () at (.5,1) {$\be$};
\end{scope}

\begin{scope}[shift={(0,-4)}]	
\colrect17{2.5}5{green}
\colrecthigh{7.5}{10.5}{green}
\colrect{11}{18}{12}{16}{green}
\colrectlow12{green}
\colrectlow{5.5}{6.5}{green}
\colrectlow7{7.5}{green}
\colrectlow8{8.5}{green}
\colrectlow9{10}{green}
\colrectlow{10.5}{11.5}{green}
\colrectlow{16.5}{18}{green}
\draw[|-|] (.5,2)--(.5,0);
\node[left] () at (.5,1) {$\al$};
\end{scope}

\end{tikzpicture}
\caption{Geometrical interpretation of $\al\leq_r\be$.  Here $\al=\id_\ve\be$, where $\ve\in\Eq(X)$ is such that each $\ve$-class is a union of $\ker(\be)$-classes.  See Remark \ref{r:leqr} for more details.}
\label{f:leqr}
\end{center}
\end{figure}
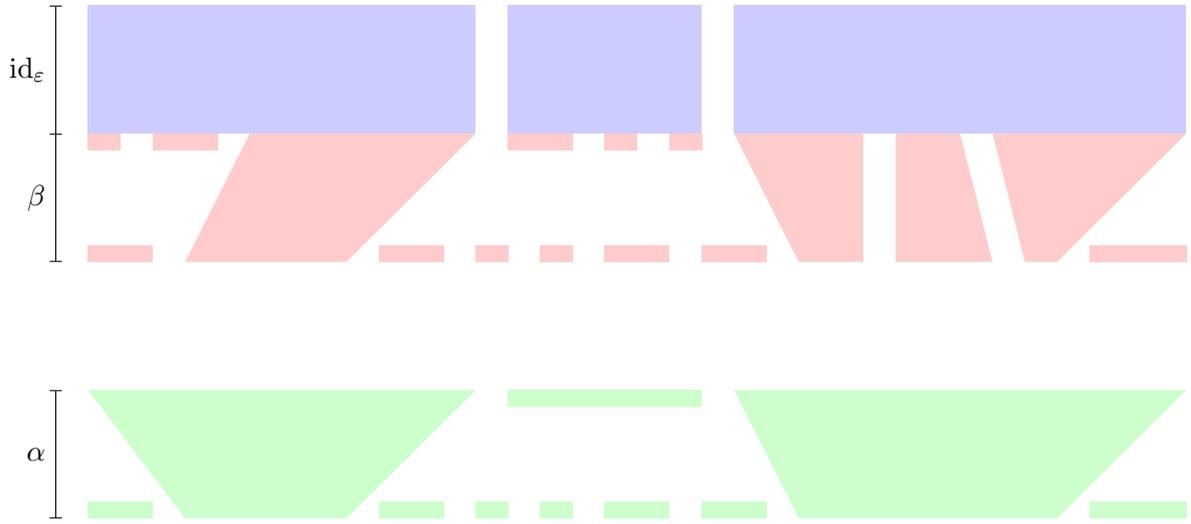

\begin{rem}
Since \ref{R1} and \ref{L1} hold for the semilattice $E=E(\I_X)$, the relations $\leq_r'$ and~$\leq_l'$ defined by
\[
\al\leq_r'\be \iff \al\in E\be \AND \al\leq_l'\be \iff \al\in \be E
\]
are still partial orders on $\P_X$.  (This also follows simply from $E$ being a semilattice containing the identity.)  It is also possible to characterise these partial orders, though the description is not as neat as those for $\leq_r$ and $\leq_l$ given in Proposition \ref{p:leqr}.

To describe the $\leq_r'$ relation ($\leq_l'$ is dual), consider partitions $\al,\be\in\P_X$, and write $\be = \begin{partn}{2} A_i&C_j\\ \hhline{~|-} B_i&D_k\end{partn}$.  Then $\al\leq_r'\be\iff\al=\id_H\be$ for some $H\sub X$.  Analysing the product $\id_H\be$, we can describe the blocks of $\al$:
\ben
\item \label{leq'1} Each lower non-transversal $D_k'$ of $\be$ is a block of $\al$.
\item \label{leq'2} For each upper non-transversal $C_j$ of $\be$, write $H_j=H\cap C_j$.  If $H_j$ is non-empty, then it is a block of $\id_H\be=\al$.  For any $x\in C_j\sm H$, $\al$ also contains the singleton $\{x\}$.
\item \label{leq'3} For each transversal $A_i\cup B_i'$ of $\be$, write $H_i=H\cap A_i$.  If $H_i$ is non-empty, then $H_i\cup B_i'$ is a block of $\al$; otherwise, $B_i'$ is a block of $\al$.  For any $x\in A_i\sm H$, $\al$ also contains the singleton~$\{x\}$.
\een
In simpler terms, this all says that any $\al\leq_r'\be$ is formed by removing any collection of upper vertices from blocks of $\be$, and replacing them by (upper) singletons.


It is also easy to see that the converse of the above is true: If conditions \ref{leq'1}--\ref{leq'3} hold, then $\al\leq_r'\be$, as $\al=\id_H\be$ for $H=\bigcup_{i\in I}H_i\cup\bigcup_{j\in J}H_j$.  Note also that conditions \ref{leq'1}--\ref{leq'3} imply $\al\pre\be$, in contrast to $\al\leq_r\be\implies\be\pre\al$ (cf.~Proposition \ref{p:leqr}).
\end{rem}

\subsection{Substructures}

Since $\P_X$ is $F$-Ehresmann, we wish to identify the monoid $\HF$-classes (cf.~Lemma~\ref{l:HeE}), and the inverse submonoid $\Reg_F(\P_X)$ consisting of all $F$-regular elements (cf.~Lemma~\ref{l:Reg}).

\begin{prop}\label{p:subs}
With respect to the semilattice $F=E(\J_X)\sub\P_X$:
\ben
\item \label{i:subs1} The monoid $H_{\id_\ve}^F$, $\ve\in\Eq(X)$, is isomorphic to the symmetric inverse monoid $\I_{X/\ve}$.
\item \label{i:subs2} The inverse monoid $\Reg_F(\P_X)$ is the dual symmetric inverse monoid $\J_X$.
\een
\end{prop}

\pf
\firstpfitem{\ref{i:subs1}} Write $H = H_{\id_\ve}^F = \set{\al\in\P_X}{\ker(\al)=\coker(\al)=\ve}$; cf.~Proposition \ref{p:RELEPX}.
As in Remark \ref{r:part_bb}, the elements of $H$ can be identified with partial bijections $X/\ve\to X/\ve$: i.e., with the elements of $\I_{X/\ve}$.  It is easy to see that multiplication in $H$ matches that in $\I_{X/\ve}$.

\pfitem{\ref{i:subs2}}  To be $\R$-related to $\id_\ve\in F$ is to have the same kernel and domain as $\id_\ve$ (cf.~Proposition~\ref{p:RL}\ref{i:RL3}), and these are $\ve$ and $X$, respectively.  Thus, to be $\R$-related to \emph{some} element of $F$ is to have domain $X$.  Similarly, to be $\L$-related to some element of $F$ is to have codomain $X$.  Thus, $\Reg_F(\B_X) = \set{\al\in\P_X}{\dom(\al)=\codom(\al)=X} = \J_X$.
\epf

The $\RF$- and $\LF$-classes of the identity element $\id_\De=\id_X$ are the monoids
\[
R_{\id_\De}^F = \set{\al\in\P_X}{\ker(\al)=\De} \AND L_{\id_\De}^F = \set{\al\in\P_X}{\coker(\al)=\De},
\]
respectively.  The $\HF$-class of $\id_\De$ is of course the symmetric inverse monoid $\I_X$; cf.~Proposition~\ref{p:subs}\ref{i:subs1}.

\begin{rem}\label{r:E}
Recall that the semilattice $E=E(\I_X)$ does not give rise to an $E$-Ehresmann structure on $\P_X$; cf.~Remark \ref{r:notE}.  For $A\in\Sub(X)$, the $\HE$-class
\[
H_{\id_A}^E = \set{\al\in\P_X}{\supp(\al)=\cosupp(\al)=A}
\]
is not a submonoid for $|A|\geq3$.  For example, consider the partitions $\al = \custpartn{1,2,3}{1,2,3}{\stline11\stline12\stline23\stline33\uarc23 \darc12}$ and $\be=\custpartn{1,2,3}{1,2,3}{\stline33\uarc12 \darc12}$, both from $\P_3$.  Then $\al$ and $\be$ both belong to the $\HE$-class of the identity, yet $\al\be=\custpartn{1,2,3}{1,2,3}{\stline23\stline33\uarc23 \darc12}$ does not.

On the other hand, $\Reg_E(\P_X) = {\set{\al\in\P_X}{\ker(\al)=\coker(\al)=\De}}$ is the symmetric inverse monoid $\I_X$.  In particular $\Reg_E(\P_X)$ is still an inverse submonoid, even though $\P_X$ is not $E$-Ehresmann.  

The $\RE$- and $\LE$-classes of $\id_X$ are not submonoids for $|X|\geq3$.
\end{rem}

\subsection{Restriction subsemigroups}\label{ss:RS}

The next lemma will help us identify the largest left- and right-restriction subsemigroups of~$\P_X$, with respect to the semilattice $F=E(\J_X)$.  For the statement and later use, we write $\nab = X\times X$ for the universal relation on $X$.  The proof utilises the idempotent $\id_\nab\in F$, which has a single block,~$X\cup X'$.  

\begin{lemma}\label{l:RestPX}
For $\al\in\P_X$, we have
\ben
\item \label{i:RestPX1} $F\al\sub \al F \iff [\dom(\al)=X$ or $\ker(\al)=\nab]$,
\item \label{i:RestPX2} $\al F\sub F\al \iff [\codom(\al)=X$ or $\coker(\al)=\nab]$.
\een
\end{lemma}

\pf
As usual, it suffices to prove \ref{i:RestPX1}.

\pfitem{($\Rightarrow$)}  Aiming to prove the contrapositive, suppose $\dom(\al)\not=X$ and $\ker(\al)\not=\nab$.  Together these imply that $\al$ has some upper non-transversal $A\subsetneq X$.  First note that ${\ker(\id_\nab\al)=\nab}$.  On the other hand, any element of $\al F$ contains the block $A$, so has non-universal kernel.  Thus, $\id_\nab\al \in F\al\sm\al F$.

\pfitem{($\Leftarrow$)}  If $\ker(\al)=\nab$, then $F\al = \{\al\} \sub \al F$.  For the rest of the proof, we assume that ${\dom(\al)=X}$.  We may then write $\al=\begin{partn}{2} A_i&\\ \hhline{~|-} B_i&C_j\end{partn}$.  (This use of the tabular notation indicates that~$\al$ has no upper non-transversals.)  Let $\ve\in\Eq(X)$, and write $X/\ve=\set{D_k}{k\in K}$.  The proof will be complete if we can show that $\id_\ve\al = \al\id_\eta$ for some $\eta\in\Eq(X)$.  As in the proof of Proposition \ref{p:leqr}, we have $\id_\ve\al = \id_{\ve\vee\ker(\al)}\al$, so we may assume that every $\ve$-class is a union of $\ker(\al)$-classes.  Thus, there is a partition $\set{I_k}{k\in K}$ of $I$ such that $D_k = \bigcup_{i\in I_k}A_i$ for all $k$.  But then $\id_\ve\al=\al\id_\eta$, where $\eta\in\Eq(X)$ has classes $C_j$ for each $j\in J$, and $\bigcup_{i\in I_k}B_i$ for each $k\in K$.  This last step is shown schematically in Figure \ref{f:RestPX}.
\epf

\begin{figure}[h]
\begin{center}
\begin{tikzpicture}[scale=.85]

\begin{scope}[shift={(0,2)}]	
\colrect1717{blue}
\colrect{7.5}{10.5}{7.5}{10.5}{blue}
\colrect{11}{18}{11}{18}{blue}
\draw[|-] (.5,2)--(.5,0);
\node[left] () at (.5,1) {$\id_\ve$};
\end{scope}

\begin{scope}[shift={(0,0)}]	
\colrect13{2.5}5{red}
\colrect{3.5}7{5.5}{6.5}{red}
\colrect{7.5}{10.5}9{10}{red}
\colrect{11}{13}{12}{13}{red}
\colrect{13.5}{14.5}{13.5}{15}{red}
\colrect{15}{18}{15.5}{16}{red}
\colrectlow12{red}
\colrectlow7{7.5}{red}
\colrectlow8{8.5}{red}
\colrectlow{10.5}{11.5}{red}
\colrectlow{16.5}{18}{red}
\draw[|-|] (.5,2)--(.5,0);
\node[left] () at (.5,1) {$\al$};
\end{scope}

\begin{scope}[shift={(0,-4)}]	
\colrect13{2.5}5{red}
\colrect{3.5}7{5.5}{6.5}{red}
\colrect{7.5}{10.5}9{10}{red}
\colrect{11}{13}{12}{13}{red}
\colrect{13.5}{14.5}{13.5}{15}{red}
\colrect{15}{18}{15.5}{16}{red}
\colrectlow12{red}
\colrectlow7{7.5}{red}
\colrectlow8{8.5}{red}
\colrectlow{10.5}{11.5}{red}
\colrectlow{16.5}{18}{red}
\draw[|-|] (.5,2)--(.5,0);
\node[left] () at (.5,1) {$\al$};
\end{scope}

\begin{scope}[shift={(0,-6)}]	
\foreach \x/\y in {1/2,2.5/6.5,7/7.5,8/8.5,9/10,10.5/11.5,12/16,16.5/18} {\colrect\x\y\x\y{green}}
\draw[-|] (.5,2)--(.5,0);
\node[left] () at (.5,1) {$\id_\eta$};
\end{scope}

\end{tikzpicture}
\caption{Verification that $\id_\ve\al=\al\id_\eta$ from the proof of Lemma \ref{l:RestPX}.}
\label{f:RestPX}
\end{center}
\end{figure}

In the next statement, we write $\ze=\{X,X'\}$ for the (unique) partition from $\P_X$ with ${\dom(\ze)=\codom(\ze)=\emptyset}$ and $\ker(\ze)=\coker(\ze)=\nab$.  We write $\sqcup$ for disjoint union.

\begin{prop}\label{p:RestPX}
With respect to the semilattice $F=E(\J_X)$:
\ben
\item \label{i:pRestPX1} $\Rest_R(\P_X,F) = \set{\al\in\P_X}{\dom(\al)=X \text{ or } \ker(\al)=\nab}$
\item[] $\phantom{\Rest_R(\P_X,F)} = \set{\al\in\P_X}{\dom(\al)=X} \sqcup \set{\al\in\P_X}{\dom(\al)=\emptyset,\ \ker(\al)=\nab}$,
\item \label{i:pRestPX2} $\Rest_L(\P_X,F) = \set{\al\in\P_X}{\codom(\al)=X \text{ or } \coker(\al)=\nab}$
\item[] $\phantom{\Rest_L(\P_X,F)} = \set{\al\in\P_X}{\codom(\al)=X} \sqcup \set{\al\in\P_X}{\codom(\al)=\emptyset,\ \coker(\al)=\nab}$,
\item \label{i:pRestPX3} $\Rest(\P_X,F) = \J_X\cup\{\ze\}$ is the dual symmetric inverse monoid $\J_X$ with an additional zero~$\ze$ attached.
\een
\end{prop}

\pf
Part \ref{i:pRestPX1} follows from Lemma \ref{l:RestPX}\ref{i:RestPX1}, in light of the fact that
\[
\big[\dom(\al)=X \text{ or } \ker(\al)=\nab\big] \iff \big[\dom(\al)=X \text{ or } [\dom(\al)=\emptyset \text{ and } \ker(\al)=\nab]\big].
\]
(For the forwards implication, observe that if $\al\in\P_X$ satisfies $\dom(\al)\not=X$ and $\ker(\al)=\nab$, then in fact $\dom(\al)=\emptyset$.)
Part \ref{i:pRestPX2} is dual to \ref{i:pRestPX1}, and then \ref{i:pRestPX3} follows from \ref{i:pRestPX1} and \ref{i:pRestPX2}.
\epf

Since we now wish to discuss the subsemigroups $\Rest_R(\P_X,F)$ and $\Rest_L(\P_X,F)$, it will be convenient to denote these by $\RR_X$ and $\LL_X$, respectively.  In what follows, we focus exclusively on $\RR_X$, as statements for $\LL_X=\RR_X^\circ$ are dual.

Consider the sets
\[
\Pfd_X = \set{\al\in\P_X}{\dom(\al)=X} \AND \Pfk_X = \set{\al\in\P_X}{\ker(\al)=\nab}
\]
of all ``full-domain'' and ``full-kernel'' partitions from $\P_X$, and also the sets
\[
D_1 = \set{\al\in\P_X}{\dom(\al)=X,\ \ker(\al)=\nab} \ANd D_0 = \set{\al\in\P_X}{\dom(\al)=\emptyset,\ \ker(\al)=\nab}.
\]
Proposition \ref{p:RestPX}\ref{i:pRestPX1} says that $\RR_X = \Pfd_X\cup\Pfk_X = \Pfd_X\sqcup D_0$, and it follows from the parenthetical observation in the proof that $\Pfk_X = D_0\sqcup D_1$.  The following structural properties are easily checked (a \emph{right zero} of a semigroup~$S$ is an element $z\in S$ such that $az=z$ for all $a\in S$):
\bit
\item $\Pfd_X$, $\Pfk_X$, $D_0$ and $D_1$ are all subsemigroups of $\P_X$, with $\Pfd_X$ a submonoid.
\item $D_1$ is the (unique) minimal ideal of $\Pfd_X$, and each element of $D_1$ is a right zero of $\Pfd_X$.
\item $D_0$ is the (unique) minimal ideal of $\RR_X$, and each element of $D_0$ is a right zero of $\RR_X$.
\eit

Further structural information concerning the monoids $\RR_X$ and $\Pfd_X$ can be found by studying their Green's relations, as defined (for any semigroup) in Subsection \ref{ss:E}.  For the next statement, recall that a semigroup $S$ is \emph{regular} if for every $x\in S$ we have $x=xax$ for some $a\in S$.  

\begin{prop}\label{p:RX_reg}
The monoids $\RR_X$ and $\Pfd_X$ are both regular.
\end{prop}

\pf
Since $\RR_X = \Pfd_X\cup D_0$, and since every element of $D_0$ is idempotent (and hence regular), it suffices to show that $\Pfd_X$ is regular.  To do so, let $\al\in\Pfd_X$ be arbitrary; we must show that $\al=\al\be\al$ for some $\be\in\Pfd_X$.  Write $\al=\begin{partn}{2} A_i&\\ \hhline{~|-} B_i&C_j\end{partn}_{i\in I,\ j\in J}$.  Since $\dom(\al)=X$ we have $I\not=\emptyset$, so we may fix some $q\in I$.  Also put $C=\bigcup_{j\in J}C_j$, but note that $C$ may be empty.  It is then easy to check that $\al=\al\be\al$ for $\be = \begin{partn}{3} B_q\cup C&B_i & \ \ \ \ \\ \hhline{~|~|-} A_q&A_i& \end{partn}_{i\in I\sm\{q\}}$, and we note that $\be$ belongs to $\Pfd_X$.  Indeed, examining the blocks of $\al$, we have $X=\bigcup_{i\in I}B_i\cup\bigcup_{j\in J}C_j = \bigcup_{i\in I}B_i\cup C$.
\epf

We may now describe Green's relations on $\RR_X$ and $\Pfd_X$.  

\begin{prop}\label{p:RX_Green}
Let $M$ be either of the monoids $\RR_X$ or $\Pfd_X$, and let $\al,\be\in M$.  Then
\ben
\item \label{i:RX1} $\al\R\be \iff [\dom(\al)=\dom(\be)$ and $\ker(\al)=\ker(\be)]$,
\item \label{i:RX2} $\al\L\be \iff [\codom(\al)=\codom(\be)$ and $\coker(\al)=\coker(\be)]$,
\item \label{i:RX3} $\al\D\be \iff \al\JJ\be \iff \rank(\al)=\rank(\be)$.
\een
\end{prop}

\pf
\firstpfitem{\ref{i:RX1} and \ref{i:RX2}}  These follow from Propositions \ref{p:RL} and \ref{p:RX_reg}.  (It is a standard result that if~$S$ is a regular subsemigroup of a semigroup $T$, then the $\R$ and $\L$ relations on $S$ are precisely the restrictions to $S$ of the $\R$ and $\L$ relations on $T$; see \cite[Proposition~2]{AHK1965}.)

\pfitem{\ref{i:RX3}}  Since ${\D}\sub{\JJ}$ in any semigroup, we have $\al\D\be\implies\al\JJ\be$.  If $\al\JJ\be$ in $M$, then $\al\JJ\be$ in the full partition monoid $\P_X$, and so $\rank(\al)=\rank(\be)$ by Proposition \ref{p:RL}\ref{i:RL6}.  It remains to show that $\rank(\al)=\rank(\be) \implies \al\D\be$ (in $M$).  So suppose $\rank(\al)=\rank(\be)$.  If $\rank(\al)=\rank(\be)=0$ (this case only arises if $M=\RR_X$), then $\dom(\al)=\dom(\be)=\emptyset$ and $\ker(\al)=\ker(\be)=\nab$, so that $\al\R\be$ by part \ref{i:RX1}, and certainly $\al\D\be$.  So suppose instead that $\rank(\al)=\rank(\be)\geq1$, and write $\al = \begin{partn}{2} A_i&\\ \hhline{~|-} B_i&C_j\end{partn}$ and $\be = \begin{partn}{2} F_i&\\ \hhline{~|-} G_i&H_k\end{partn}$.  Then using parts \ref{i:RX1} and \ref{i:RX2}, we have $\al\R\ga\L\be$ for $\ga = \begin{partn}{2} A_i&\\ \hhline{~|-} G_i&H_k\end{partn}\in\Pfd_X\sub M$, and so $\al\D\be$.
\epf

\begin{rem}
Proposition \ref{p:RX_Green}\ref{i:RX1} of course simplifies for $M=\Pfd_X$, given that all partitions in this monoid have domain $X$.  In $\Pfd_X$ we have $\al\R\be\iff\ker(\al)=\ker(\be)$.
\end{rem}

In fact, it follows from regularity of $\RR_X$ and $\Pfd_X$ (and \cite[Lemma 3]{EH2020}) that Green's $\leqR$ and~$\leqL$ preorders on $\RR_X$ and $\Pfd_X$ are as described in Proposition \ref{p:RL}.  Among other things, the next result shows that this is also the case for $\leqJ$.  For the statement, and the following discussion, it is more convenient to speak of the partial order on $\JJ$-classes.  (The $\JJ$-classes of an arbitrary semigroup $S$ are partially ordered, for $x,y\in S$, by $J_x\leq J_y \iff x\leqJ y$.)  The notation for the $\D$-classes in the statement explains our earlier names for the sets $D_0$ and~$D_1$.

\begin{prop}\label{p:RX_Green2}
Let $M$ be either of the monoids $\RR_X$ or $\Pfd_X$, and let
\[
z = \begin{cases}
0 &\text{if $M=\RR_X$}\\
1 &\text{if $M=\Pfd_X$.}
\end{cases}
\]
\ben
\item \label{i:RX4}  The ${\D}={\JJ}$-classes of $M$ are the sets 
\[
D_\mu = \set{\al\in M}{\rank(\al)=\mu} \qquad\text{for each cardinal $z\leq\mu\leq |X|$,}
\]
and these form a chain under the usual partial order on $\JJ$-classes: $D_\mu\leq D_\nu \iff \mu\leq\nu$.
\item \label{i:RX5}  For any cardinal $z\leq\mu\leq|X|$, each group $\H$-class contained in $D_\mu$ is isomorphic to the symmetric group $\S_\mu$.
\een
\end{prop}

\pf
\firstpfitem{\ref{i:RX4}}  The claim concerning the $\D$-classes follows from Proposition \ref{p:RX_Green}\ref{i:RX3}.  For the second claim, let $\al\in D_\mu$ and $\be\in D_\nu$; we must show that $\al\leqJ\be \iff \mu\leq\nu$.  For the forwards implication, we have 
\[
\al\leqJ\be \text{ (in $M$)} \implies \al\leqJ\be \text{ (in $\P_X$)} \implies \mu\leq\nu,
\]
using Proposition \ref{p:RL}\ref{i:RL5} in the last step.  For the converse, suppose $\mu\leq\nu$.  If $\mu=0$, then $\al=\be\al$ and so $\al\leqJ\be$.  So now suppose $\mu\geq1$, and write $\al = \begin{partn}{2} A_i&\\ \hhline{~|-} B_i&C_j\end{partn}$ and $\be = \begin{partn}{2} F_k&\\ \hhline{~|-} G_k&H_l\end{partn}$.
Since $|I|=\mu\leq\nu=|K|$, we may assume without loss of generality that $I\sub K$.  Since $\mu\geq1$, we may fix some $q\in I$, and we also write $H=\bigcup_{l\in L}H_l\cup\bigcup_{k\in K\sm I}G_k$.  Then we have $\al=\ga\be\de$ for $\ga = \begin{partn}{2} A_i& \\ \hhline{~|-} F_i&F_k\end{partn}_{i\in I,\ k\in K\sm I}$ and $\de = \begin{partn}{3} G_q\cup H&G_i&\\ \hhline{~|~|-} B_q&B_i&C_j\end{partn}_{i\in I\sm\{q\}}$,
both from $\Pfd_X$.  This shows that $\al\leqJ\be$, and completes the proof.

\pfitem{\ref{i:RX5}}  Comparing Propositions \ref{p:RL} and \ref{p:RX_Green}, the $\H$-class in $M$ of any idempotent $\al\in D_\mu$ coincides with its $\H$-class in $\P_X$.  By \cite[Lemma 2.3]{ERPX}, the latter is isomorphic to $\S_\mu$.
\epf

\begin{rem}
Since the ideals of a semigroup $S$ are in one-one correspondence with the downward-closed subsets of the poset $(S/{\JJ},\leq)$, it follows from Proposition \ref{p:RX_Green2}\ref{i:RX4} that the ideals of the monoids $\RR_X$ and $\Pfd_X$ form chains.
\end{rem}

\begin{rem}
Figure \ref{f:RP4} shows an \emph{egg-box diagram} of $\RR_4 = \Rest_R(\P_4,F)$, produced using the {\sc Semigroups} package for GAP \cite{GAP,GAP4}.  In such a diagram, $\R$- or $\L$-related elements are placed in the same row or column, respectively, so that cells are $\H$-classes.  The large boxes are the ${\D}={\JJ}$-classes, ordered $D_0<D_1<\cdots<D_4$ (bottom to top).
The group $\H$-classes are shaded (grey or dark orange) in Figure \ref{f:RP4}, and labelled with their standard GAP names.  The elements of $\J_4$ are shaded orange, so the dark orange $\H$-classes contain the elements of the semilattice $F=E(\J_4)$.  The green $\H$-class contains $\ze$, defined above.

The minimal ideal $D_0$ of~$\RR_4$ is the bottom $\D$-class in Figure~\ref{f:RP4}.  An egg-box diagram of~${\Pfd_4 = \RR_4\sm D_0}$ may be obtained by removing the bottom $\D$-class; the minimal ideal~$D_1$ of~$\Pfd_4$ is then the bottom $\D$-class remaining.  The bottom two $\D$-classes alone constitute the semigroup~${\Pfk_4=D_0\cup D_1}$.
\end{rem}

\begin{figure}[h]
\begin{center}
\includegraphics[width=\textwidth]{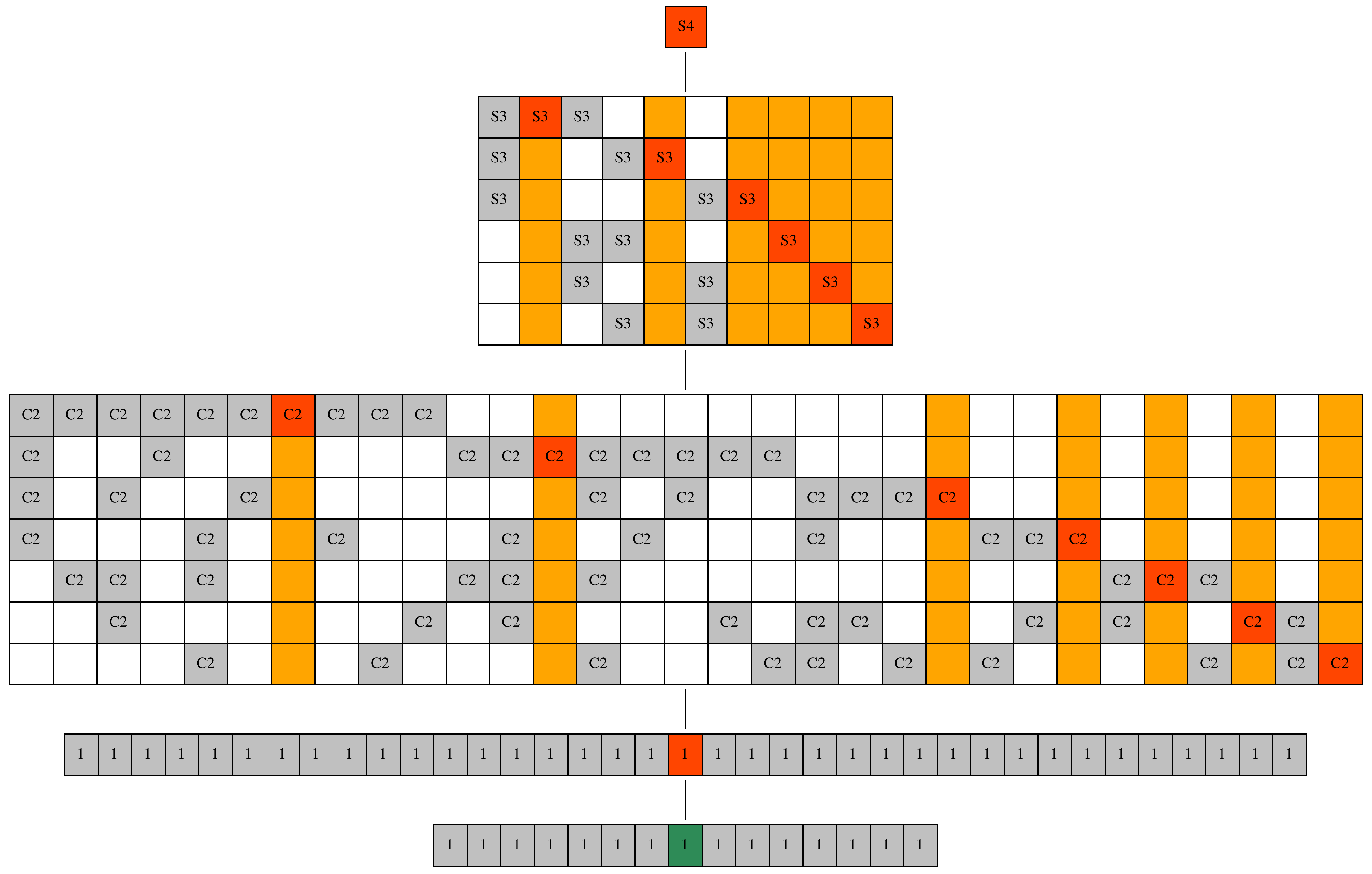}
\caption{Egg-box diagram of the largest right-restriction subsemigroup $\RR_4 = \Rest_R(\P_4,F)$ of the partition monoid~$\P_4$, with respect to the semilattice $F=E(\J_4)$.  The elements of $\J_4$ are shaded orange; the dark orange $\H$-classes contain the elements of $F$; the green $\H$-class contains~$\ze$.  For more details see Subsection \ref{ss:RS}.}
\label{f:RP4}
\end{center}
\end{figure}

Proposition \ref{p:RestPX} describes the largest left-, right- and two-sided restriction subsemigroups of~$\P_X$.  Because the monoid $\Pfd_X$ and its dual are so natural, it seems worth recording the following.

\begin{prop}\label{p:Pfd}
The monoids
\[
\Pfd_X = \set{\al\in\P_X}{\dom(\al)=X} \AND \Pfcd_X = \set{\al\in\P_X}{\codom(\al)=X} 
\]
of all full-domain and full-codomain partitions are right- and left-restriction monoids, respectively, with respect to the semilattice $F=E(\J_X)$ coming from the dual symmetric inverse monoid over $X$.  Further, $\Pfd_X\cap\Pfcd_X = \J_X$.  \epfres
\end{prop}

\begin{rem}\label{r:Pfd}
The monoid $\Pfd_X$ was independently discovered by Stokes \cite{Stokes_prep} in his work on categories and constellations.  We believe this monoid is of significant independant interest, and worthy of further study, but this is beyond the scope of the current work.  However, we make one final comment on the duality between $\Pfd_X$ and $\PT_X$.  Namely, we have factorisations
\[
\PT_X = \la \T_X\cup\I_X\ra = E(\I_X)\cdot\T_X \AND \Pfd_X = \la\T_X\cup\J_X\ra = \T_X \cdot E(\J_X).
\]
Thus, while elements of $\PT_X$ can be thought of full transformations (from $\T_X$) ``restricted in the domain'' by subsets (elements of $E(\I_X)$), elements of $\Pfd_X$ can be thought of as full transformations ``restricted in the codomain'' by quotients (elements of $E(\J_X)$).
\end{rem}


We conclude this subsection by considering the monoid $\RR_X$ as an $F$-Ehresmann monoid in its own right.  First, the $\HF$-class in $\RR_X$ of $\id_\ve$, $\ve\in\Eq(X)$, is the monoid
\[
H_{\id_\ve}^F = \set{\al\in\RR_X}{\ker(\al)=\coker(\al)=\ve}.
\]
When $\ve=\nab$ is the universal equivalence, $H_{\id_\nab}^F = \{\id_\nab,\ze\}$, where $\ze=\{X,X'\}$ was defined above.  Otherwise,
\[
H_{\id_\ve}^F = \set{\al\in\P_X}{\dom(\al)=X,\ \ker(\al)=\coker(\al)=\ve} \qquad\text{for $\ve\not=\nab$.}
\]
As in Remark \ref{r:part_bb}, this is isomorphic to the monoid of injective (total) functions $X/\ve\to X/\ve$.  When $X/\ve$ is finite, this is of course the symmetric group $\S_{X/\ve}$, and here $H_{\id_\ve}^F=H_{\id_\ve}$ is precisely the usual Green's $\H$-class of $\id_\ve$.  When $|X|\geq2$, the $\RF$- and $\LF$-classes in~$\RR_X$ of the identity element $\id_\De=\id_X$ are
\begin{align*}
R_{\id_\De}^F &= \set{\al\in\RR_X}{\ker(\al)=\De} = \set{\al\in\P_X}{\dom(\al)=X,\ \ker(\al)=\De} \\[2mm]
\text{and}\qquad L_{\id_\De}^F &= \set{\al\in\RR_X}{\coker(\al)=\De} = \set{\al\in\P_X}{\dom(\al)=X,\ \coker(\al)=\De}.
\end{align*}
The second of these is the full transformation semigroup $\T_X$ (cf.~Subsection \ref{ss:Green}), while the first is the $\R$-class in $\P_X$ of the identity element (cf.~Proposition \ref{p:RL}\ref{i:RL3}): i.e., the monoid of all right-invertible partitions from $\P_X$.  (This latter submonoid played a key role in \cite{JEipms}, where it was denoted $\LL_X$, conflicting with our notation above.  Indeed, continuing just for now with the~$\LL_X$ notation of \cite{JEipms}, and writing $\RR_X=\LL_X^\circ=\set{\al\in\P_X}{\codom(\al)=X,\ \coker(\al)=\De}$, we have $\P_X=\LL_X\RR_X$ for infinite $X$; cf.~\cite[Remark 11]{JEipms}.  When $X$ is finite, $\LL_X=\RR_X=\S_X$.)

The discussion of the previous paragraph applies to the semigroup $\Pfd_X$ as well, but without the need to single out the $\ve=\nab$ case, as the $\D$-class $D_0$ containing $\ze$ has been removed.
It is again easy to see that $\Reg_F(\RR_X)=\Reg_F(\Pfd_X)=\J_X$.

We conclude this subsection with a list of open questions, all of which beyond the scope of the current work:
\bit
\item Describe the idempotent-generated subsemigroup of $\Pfd_X$.  Does this coincide with the singular ideal in the finite case?  Cf.~\cite{Howie1966,EF2012,East2011_2,Howie1978}.
\item In the finite case, are the proper ideals of $\Pfd_X$ idempotent-generated?  What are the ranks (and idempotent ranks, if appropriate)?  Cf.~\cite{HM1990,EG2017,Gray2008,Gray2014}.
\item Describe the biordered set of idempotents of $\Pfd_X$, and the structure of the free semigroup over this biorder.  Cf.~\cite{GR2012,GR2012_2,DG2014,DGR2017,Nambooripad1979,Easdown1985}.
\item Give a presentation for $\Pfd_X$ in the finite case.  Cf.~\cite{Aizenstat1958,East2018,East2011}.
\item Determine the congruences of $\Pfd_X$.  Cf.~\cite{EMRT2018,Malcev1952,ERPX,Klimov1977,ER2020}.
\item Does every (regular) right-restriction Ehresmann semigroup embed as a $(2,1,1)$-subalgebra of some $\Pfd_X$?  Cf.~\cite{MS2020,Gould_notes}.  (The same question can be asked for other signatures.)
\item Describe the simple/projective/injective $\Pfd_X$-modules.  Cf.~\cite{MS2020,Stein2017,Stein2016}.
\eit

\subsection{Categories and representations}\label{ss:CRP}

Let us now return our attention to the partition monoid $\P_X$, considered as an $F$-Ehresmann semigroup with respect to the semilattice $F=E(\J_X)=\bigset{\id_\ve}{\ve\in\Eq(X)}$ of the dual symmetric inverse monoid~$\J_X$.  In particular, we will investigate $C=\bC(\P_X,F)$, the (Ehresmann) category whose construction was detailed in Subsection \ref{ss:cats}.  As with the monoid $\B_X$ of binary relations (cf.~Subsection \ref{ss:catBX}), this category has a neat combinatorial description.

Identifying an idempotent $\id_\ve$ from $F$ with the partition $X/\ve$, we will identify the object set of the category $C=\bC(\P_X,F)$ with $\Part(X)$, the set of all partitions of $X$.  ($\Part(X)$ is not to be confused with~$\P_X$ itself, which consists of partitions of $X\cup X'$.)  As ever, the morphism set
\[
C(\bA,\bB) = \set{\al\in\P_X}{X/\ker(\al)=\bA,\ X/\coker(\al)=\bB} \qquad\text{for $\bA,\bB\in\Part(X)$}
\]
is (identified with) the set of all partial bijections $\bA\to\bB$.

The category $C'=\bC(\RR_X,F)$ associated to the right-restriction monoid $\RR_X=\Rest_R(\P_X,F)$ also has object set $\Part(X)$, but now the morphism set
\begin{align*}
C'(\bA,\bB) &= \set{\al\in\RR_X}{X/\ker(\al)=\bA,\ X/\coker(\al)=\bB}
\end{align*}
consists of all injective (total) functions $\bA\to\bB$, apart from when $\bA$ is the trivial partition $\bX = X/\nab=\{X\}$; in this latter case, $C'(\bX,\bB)$ consists of all total functions $\bX\to\bB$ (which are of course injective by default) and also the unique empty map $\bX\to\bB$.  

The category $C''=\bC(\Pfd_X,F)$ associated to the full-domain partition monoid $\Pfd_X$ has essentially the same description as $C'$, except that no extra trouble is required for morphisms $\bX\to\bB$; these are only the (injective) total functions.

When $X$ is finite, the $\leq_l$ partial order on the right-restriction monoids $\RR_X$ and $\Pfd_X$ is of course finite-below, so Theorem \ref{t:Stein1} links the representation theory of these monoids with that of the categories $C'$ and $C''$.  We believe it would be interesting to investigate this, perhaps using the FI-modules of Church, Ellenberg and Farb \cite{CEF2015,CE2017}, but this is beyond the scope of the current work.  Note for example that $C''$ is equivalent (in the categorical sense) to the category whose objects are the non-empty subsets of $X$, and whose morphisms are injections.

Finally, for finite $\bA\in\Part(X)$, the endomorphism monoid $C''(\bA,\bA)$ is precisely the symmetric group $\S_{\bA}$.  It follows that $C''=\bC(\Pfd_X,F)$ is an EI-category when $X$ is finite.  Theorem~\ref{t:Stein2} (and Proposition \ref{p:Pfd}) then tells us that for any field $\K$ of suitable characteristic, the maximal semisimple image of $\K[\Pfd_X]$ is precisely $\K[\J_X]$.  Given that the group $\H$-classes of $\Pfd_X$ are symmetric groups $\S_n$ for $n\leq|X|$, the restriction on the characteristic of $\K$ is that it should be either~$0$ or else greater than $|X|$.

The previous paragraph concerned the category $C''=\bC(\Pfd_X,F)$.  Similar statements \emph{almost} apply to the category $C'=\bC(\RR_X,F)$, except that this category \emph{just} fails to be an EI-category.  Indeed, the endomorphism monoid $C'(\bX,\bX)$ alone contains a single non-automorphism: namely, the empty map~$\bX\to\bX$, corresponding to the partition $\ze$ defined before Proposition \ref{p:RestPX}.

\section{Other diagram monoids}\label{s:ODM}

We close with a brief discussion of two more families of diagram monoid, one contained in $\P_X$ (the partial Brauer monoid, considered in Subsection \ref{ss:PBX}), and the other containing it (the rook partition monoid, considered in Subsection \ref{ss:RPX}).

\subsection{(Partial) Brauer monoids}\label{ss:PBX}

A partition $\al\in\P_X$ is a \emph{Brauer} or \emph{partial Brauer partition} if every block $A$ of $\al$ satisfies $|A|=2$ or $|A|\leq2$, respectively.  See Figure \ref{f:B6} for some examples.  The set $\PB_X$ of all partial Brauer partitions is a submonoid of $\P_X$ (called the \emph{partial Brauer monoid}) for any $X$, while the set $\B_X$ of all Brauer partitions is a submonoid (called the \emph{Brauer monoid}) if and only if $X$ is finite \cite{JE_IBM}.  (When $X$ is infinite, $\B_X$ is a generating set for $\PB_X$, by \cite[Corollary 4.4]{JE_IBM}.)

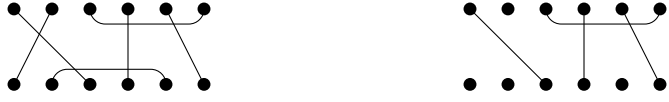
\begin{figure}[h]
\begin{center}
\begin{tikzpicture}[scale=.5]

\begin{scope}[shift={(0,0)}]	
\uvs{1,...,6}
\lvs{1,...,6}
\stline13
\stline21
\stline44
\stline56
\uarc36
\darc25
\end{scope}

\begin{scope}[shift={(12,0)}]	
\uvs{1,...,6}
\lvs{1,...,6}
\stline13
\stline44
\stline56
\uarc36
\end{scope}

\end{tikzpicture}
\caption{Sample elements of $\B_6$ (left) and $\PB_6$ (right).}
\label{f:B6}
\end{center}
\end{figure}

Since $\PB_X$ contains the symmetric inverse monoid $\I_X$, it contains in particular the semilattice $E=E(\I_X)$.  On the other hand, $\PB_X$ does not contain the semilattice $F=E(\J_X)$.  One might wonder if $\PB_X$ is $E$-Ehresmann, even though $\P_X$ itself is not.  By \cite[Theorem 4.9]{DDE2019}, Green's relations and preorders on $\PB_X$ have the same description as those on $\P_X$; cf.~Proposition~\ref{p:RL}.  It follows then that the sets $E_L(\al)$ and $E_R(\al)$, $\al\in\PB_X$, are as described in Lemma~\ref{l:eq}.  Consequently, the $\RE$ and $\LE$ relations on $\PB_X$ are as described in Proposition~\ref{p:RELEPX}.  It again quickly follows that \ref{R1} and \ref{L1} hold.  However, since the elements $\al,\be,\th$ constructed in Remark \ref{r:notE} belong to $\PB_2$, it follows that \ref{R2} and \ref{L2} still do not hold in general.  This means that $\PB_X$ is not $E$-Ehresmann either, for $|X|\geq2$.

Despite this, one may still compute the various subsets of $\PB_X$ associated to the $\RE$, $\LE$ and $\HE$ relations.  For example, the $\HE$-class of an idempotent from $E$ has the form
\[
H_{\id_A}^E = \set{\al\in\PB_X}{\supp(\al)=\cosupp(\al)=A} \qquad\text{for $A\sub X$.}
\]
This set is in one-one correspondence with the set $\B_A$ of all Brauer partitions over $A$, and is consequently a submonoid of $\PB_X$ if and only if $A$ is finite.
As in Remark \ref{r:E}, $\Reg_E(\PB_X)=\I_X$ is still an inverse submonoid.

\subsection{Rook partition monoids}\label{ss:RPX}

The \emph{rook partition monoid}, denoted $\RP_X$, is a diagram monoid containing $\P_X$.  It has been studied in a number of settings, and under a variety of different names; see for example \cite{East2018_rook,Grood2006,HR2005,Martin2000,Martin1996}.  

To give the first of two equivalent definitions of $\RP_X$, fix some symbol $\infty$ not belonging to~$X$, and let $Y=X\cup\{\infty\}$.  Then $\RP_X$ is the set of all partitions $\al$ of $Y\cup Y'$ such that $\infty$ and $\infty'$ belong to the same block of $\al$.  It is clear that $\RP_X$ is a submonoid of $\P_Y$, and in the finite case we obtain the natural tower of embeddings
\begin{equation}\label{e:PRP}
\P_0 \hookrightarrow \RP_0 \hookrightarrow \P_1 \hookrightarrow \RP_1 \hookrightarrow \P_2 \hookrightarrow \RP_2 \hookrightarrow \cdots,
\end{equation}
utilised by Halverson and Ram in their work on representations of partition algebras \cite{HR2005}; see also \cite{Martin2000,Martin1996}.

In the alternative definition, the elements of $\RP_X$ are partitions of subsets of $X\cup X'$; the elements of $X\cup X'$ ``missing'' from such a partition are those belonging to the same block as $\infty$ and $\infty'$ in the previous definition.  The equivalence of the two definitions is akin to the standard embedding of $\PT_X$ in $\T_{X\cup\{\infty\}}$ as the set of transformations mapping $\infty\mt\infty$ \cite{Vagner1956}.  In this way,~$\RP_X$ may be thought of as a ``partial version'' of $\P_X$.

A rook partition $\al\in\RP_X$ is typically drawn as a graph on vertex set $X\cup X'$, with the elements in the same block as $\infty$ and $\infty'$ indicated by white vertices; these are often called the \emph{rook dots} of $\al$.  See Figure \ref{f:RP10}, which also shows how rook partitions are multiplied: essentially as ordinary partitions, but subject to the rule that rook dots ``kill'' ordinary blocks.

\begin{figure}[H]
\begin{center}
\begin{tikzpicture}[scale=.35]
\begin{scope}[shift={(0,0)}]	
\uvws{3,9,10}
\lvws{1}
\uverts{1,2,4,5,6,7,8}
\lverts{2,3,4,5,6,7,8,9,10}
\uarcs{1/2,2/4,5/6,7/8}
\darcs{4/5,6/7,9/10}
\darcx26{.8}
\stlines{4/3,5/5,8/8}
\draw(0.5,1)node[left]{$\al=$};
\end{scope}
\begin{scope}[shift={(0,-4)}]	
\uvws{1,8,10}
\lvws{3}
\uverts{2,...,7,9}
\lverts{1,2,4,5,...,10}
\uarcs{2/3,5/6}
\darcs{1/2,4/5,6/7,8/9}
\darcx7{10}{.8}
\stlines{4/4,7/7,9/9}
\draw(0.5,1)node[left]{$\be=$};
\draw(13,3)node{$\longrightarrow$};
\end{scope}
\begin{scope}[shift={(15,-1)}]	
\uverts{1,2,4,5,6,7,8}
\lverts{2,3,4,5,6,7,8,9,10}
\uarcs{1/2,2/4,5/6,7/8}
\darcs{4/5,6/7,9/10}
\darcx26{.8}
\stlines{4/3,5/5,8/8}
\uvws{3,9,10}
\lvws{1}
\end{scope}
\begin{scope}[shift={(15,-3)}]	
\uvws{1,8,10}
\lvws{3}
\uverts{2,...,7,9}
\lverts{1,2,4,5,...,10}
\uarcs{2/3,5/6}
\darcs{1/2,4/5,6/7,8/9}
\darcx7{10}{.8}
\stlines{4/4,7/7,9/9}
\draw(13,2)node{$\longrightarrow$};
\end{scope}
\begin{scope}[shift={(30,-2)}]	
\uvws{3,7,8,9,10}
\lvws{3,8,9}
\uverts{1,2,4,5,6}
\lverts{1,2,4,5,6,7,10}
\uarcs{1/2,2/4,4/5,5/6}
\darcs{1/2,4/5,5/6,6/7,7/10}
\stlines{5/5}
\draw(10.5,1)node[right]{$=\al\be$};
\end{scope}
\end{tikzpicture}
\end{center}
\vspace{-5mm}
\caption{Multiplication of two rook partitions $\al,\be\in\RP_{10}$.}
\label{f:RP10}
\end{figure}

Note that $\RP_X$ contains the (ordinary) partition monoid $\P_X$, as the set of all rook partitions with no rook dots: i.e., the partitions of $Y\cup Y'$ containing the block $\{\infty,\infty'\}$.  In particular, it follows that~$\RP_X$ contains the two semilattices $E=E(\I_X)$ and $F=E(\J_X)$, but it turns out that $\RP_X$ is neither $E$- nor $F$-Ehressman.  This is clear for $E$, since $\RP_X$ contains $\P_X$, which is not $E$-Ehresmann.  For~$F$, consider the rook partitions
\begin{equation}\label{e:RP2}
\al = \custpartn{1}{1}{\stline11\uvws{2}\lvws{2}} \COMMA \be = \id_\De =  \custpartn{1,2}{1,2}{\stline11\stline22} \COMMA \th = \custpartn{2}{2}{\stline22\uvws{1}\lvws{1}} \COMMA \text{all from $\RP_2$.}
\end{equation}
Then $\th\al=\custpartn{}{}{\uvws{1,2}\lvws{1,2}}$ and $\th\be=\th$, and $F_L(\al)=F_L(\be)=F_L(\th\be)=\{\id_\De\}$, while $\id_\nab=\custpartn{1,2}{1,2}{\stline11\stline22\uarc12\darc12}$ belongs to $F_L(\th\al)$.  This shows that $(\al,\be)\in{\RF}$ though $(\th\al,\th\be)\not\in{\RF}$: i.e., $\RF$ is not a left congruence.

However, $\RP_X$ does have a natural Ehresmann structure, still based on dual symmetric inverse monoids.  To see this, we now return to thinking of $\RP_X$ as a submonoid of $\P_Y$ as above, with $Y=X\cup\{\infty\}$.  Note that $\RP_X$ does not contain the entire dual symmetric inverse monoid~$\J_Y$; rather, $\RP_X\cap\J_Y$ (of course) consists of all block bijections $\bA\to\bB$ from $\J_Y$ such that the block of $\bA$ containing $\infty$ is mapped to the block of $\bB$ containing $\infty$.  We will write $\RJ_X = \RP_X\cap\J_Y$ for the set of all such block bijections of $Y=X\cup\{\infty\}$, and call this the \emph{rook dual symmetric inverse monoid}.  The monoid $\RJ_X$ appears in the articles \cite{MS2019,KMU2015,KM2011}.  

Although $\RP_X$ does not contain all of $\J_Y$, it does contain $E(\RJ_Y)=E(\J_Y)$.  To avoid confusion, we will write $G=E(\J_Y)$ for this semilattice.  Since $\P_Y$ is $G$-Ehresmann (Theorem \ref{t:F}), it follows immediately that so too is~$\RP_X$.  
(With $\al,\be\in\RP_2$ as in \eqref{e:RP2}, we have $G_L(\al)=\{\al,\be\}$ and $G_L(\be)=\{\be\}$, so that $(\al,\be)\not\in{\RG}$, meaning that these rook partitions cause no problems.)

One could now calculate the various substructures and categories associated to this $G$-Ehresmann structure on $\RP_X$; for example, $\Reg_G(\RP_X)=\RJ_X$.  Of particular significance, we believe it could be fruitful to investigate the links between the various Ehresmann structures on the monoids in the tower~\eqref{e:PRP}: \emph{viz}.,
\[
\begin{array}{ccccc}
\P_X &\hookrightarrow& \RP_X &\hookrightarrow& \P_Y\\
\text{$F$-Ehresmann} && \text{$G$-Ehresmann} && \text{$G$-Ehresmann}.
\end{array}
\]
A similar approach was taken to the analogous tower of (rook) dual symmetric inverse monoids in \cite{MS2019}.  Although we hope this could lead to some insights into the representation theory of (rook) partition monoids and algebras, it is beyond the scope of the current work.

\footnotesize
\def\bibspacing{-1.1pt}
\bibliography{biblio}
\bibliographystyle{abbrv}
\end{document}